\documentclass[11pt,a4paper,leqno]{article}
\usepackage[margin=27mm]{geometry}
\usepackage[T1]{fontenc}
\usepackage{lmodern,microtype}
\usepackage{amsmath,amssymb,amsfonts,amsthm,mathtools}
\usepackage{enumitem,booktabs,array}
\usepackage[dvipsnames]{xcolor}
\usepackage[colorlinks=true,linkcolor=MidnightBlue,citecolor=MidnightBlue,urlcolor=MidnightBlue]{hyperref}
\usepackage{bookmark}

\numberwithin{equation}{section}
\newtheorem{theorem}{Theorem}[section]
\newtheorem{proposition}[theorem]{Proposition}
\newtheorem{lemma}[theorem]{Lemma}

\theoremstyle{definition}

\theoremstyle{remark}
\newtheorem{remark}[theorem]{Remark}
\newcommand{\R}{\mathbb R}
\newcommand{\Z}{\mathbb Z}
\newcommand{\E}{\mathbb E}
\newcommand{\Pp}{\mathbb P}
\newcommand{\GP}{\mathbb P^G}
\newcommand{\cE}{\mathcal E}
\newcommand{\cC}{\mathcal C}
\newcommand{\cG}{\mathcal G}
\newcommand{\cH}{\mathcal H}
\newcommand{\cL}{\mathcal L}
\newcommand{\cV}{\mathcal V}
\newcommand{\eps}{\varepsilon}
\newcommand{\capR}{\operatorname{cap}_{\R^d}}
\newcommand{\capZ}{\operatorname{cap}_{\Z^d}}
\newcommand{\dd}{\,\mathrm d}
\DeclareMathOperator{\Var}{Var}
\DeclareMathOperator{\Cov}{Cov}
\DeclareMathOperator{\supp}{supp}
\setlist[itemize]{leftmargin=*,topsep=3pt,itemsep=2pt}
\setlist[enumerate]{leftmargin=*,topsep=3pt,itemsep=3pt}
\hypersetup{pdftitle={Phase transitions in the hard-wall large deviations for the Gaussian free field}}
\title{Phase transitions in the hard-wall large deviations\\[2pt] for the Gaussian free field}
		\author{Xinyi Li\thanks{Beijing International Center for Mathematical Research, Peking University, \href{mailto:xinyili@bicmr.pku.edu.cn}{xinyili@bicmr.pku.edu.cn}}\and 
        Jialu Shi\thanks{Department of Pure Mathematics and Mathematical Statistics, University of Cambridge,  
        \href{mailto:js3007@cam.ac.uk}{js3007@cam.ac.uk}}}
\date{\today}
\begin{document}
\maketitle
\begin{abstract}
We study the probability that the Gaussian free field on $\Z^d$, $d\ge3$,
stays above a level $-h_N$ throughout the lattice $N$-blow-up  of a bounded
Lipschitz domain, when $h_N$ is close to the entropic-repulsion height
$\sqrt{4g_0\log N}$. Two mechanisms compete: a collective upward displacement
of the field, whose cost is governed by capacity, and the suppression of the
remaining local violations, whose cost is their expected number, both of which are of
order $N^{d-2}/\log N$ exactly in the window
$h_N(\beta)^2=g_0(4\log N+\log\log N+2\beta)$. We obtain a sharp asymptotic rate  expressed through a
variational problem combining a Dirichlet energy with an exponential
potential inside this critical window. Outside this window we obtain the asymptotics at the appropriate
moving speeds: capacity alone governs the cost as $\beta\to-\infty$, and the
expected number of violations as $\beta\to+\infty$. As an application, we
derive lower-deviation bounds for the cover level of a macroscopic body by random
interlacements, which determine the order of decay at the transition.
\end{abstract}

\setcounter{section}{-1}
\section{Introduction and main results}\label{sec:introduction}

Let $\varphi$ be the Gaussian free field on $\Z^d$, $d\ge3$, and let $D_N$
be the blow-up by a factor $N$ of a nice bounded domain $D\subset\R^d$. This
article concerns the probability that $\varphi$ stays above a level $-h_N$
at every site of $D_N$, when $h_N$ is close to $\sqrt{4g_0\log N}$; here
$g_0$ is the variance of $\varphi_x$. This is the height by which the field
is lifted when it is conditioned to be nonnegative on $D_N$
\cite{BDZ1995}, and it is where the mechanism behind the hard-wall event
changes.

The field has two ways of meeting the constraint. It can rise collectively,
at a Gaussian cost governed by the capacity of $D$; or it can stay
essentially where it is, in which case one pays for every site that would
otherwise fall below the wall. Our main result,
Theorem~\ref{thm:gff-window}, identifies the window of levels in which the
two costs are of the same order, $N^{d-2}/\log N$, and expresses the rate as
a variational problem combining a Dirichlet energy with an exponential
potential; these two terms are the costs of the two strategies, rising
collectively and paying for each site below the wall. We also determine the
asymptotics on either side of the window
(Theorems~\ref{thm:subcritical-sharp} and~\ref{thm:supercritical-sharp}),
and apply these estimates to the cover levels of random interlacements
(Proposition~\ref{prop:ri-phase}). Letting the wall move with $N$ (see \eqref{eq:window}) is what
makes it possible to follow the passage from a capacitary regime to a regime
governed by local violations. At the transition itself, the rate can only be
identified by keeping track of both contributions at once.
\paragraph{Background.}
The hard-wall problem models an interface that must stay above an
impenetrable substrate. To make room for its downward fluctuations, the
interface moves away from the wall as a whole; this effect is called
entropic repulsion. Lebowitz and Maes \cite{LebowitzMaes1987} obtained early
qualitative results for convex interfaces in a soft-wall setting, where the
wall is replaced by an external repulsive potential. For the infinite-volume
lattice GFF in dimension $d\ge3$, Bolthausen, Deuschel and Zeitouni
\cite{BDZ1995} determined the sharp logarithmic asymptotics of the
probability that the field is positive on a large box: it decays at speed
$N^{d-2}\log N$, with a rate given by capacity, and the conditioned field is
lifted by roughly $\sqrt{4g_0\log N}$. Deuschel and Giacomin
\cite{DeuschelGiacomin1999} then showed that, once this diverging height is
subtracted, the local law of the conditioned field converges to that of the
free field. In short, capacity describes the macroscopic displacement, and
the local fluctuations are those of the original field.

Boundary conditions matter. For the zero-boundary field on a finite box,
Deuschel \cite{Deuschel1996} showed that positivity on the whole box has a
surface-order cost $N^{d-1}$. Imposing the wall up to the boundary of a
Dirichlet box is therefore a different problem from imposing it on a large
region of the infinite-volume field, which is the setting of this article.

In two dimensions the problem is formulated for a finite-volume Dirichlet
field, with the constrained bulk kept away from the boundary. Bolthausen,
Deuschel and Giacomin \cite{BolthausenDeuschelGiacomin2001} related
entropic repulsion to the maximum of this logarithmically correlated field.
Around the usual centering of the extremes, the hard-wall probability is of
order one and is described by the extremal-process results of Biskup and
Louidor \cite{BiskupLouidor2020}, through a randomly shifted Gumbel law.
Deeper in the large-deviation regime the cost is quadratic in the height
deficit and is governed by a relative capacity; Fels, Louidor and Wu
\cite{FelsLouidorWu2026} have recently sharpened these estimates and
identified the order of the subleading correction. This planar passage from
typical extremes to collective repulsion has a window and a limiting
description different from those found here for the transient field.

Several works add disorder. Bertacchi and Giacomin
\cite{BertacchiGiacomin2002,BertacchiGiacomin2003} considered a random wall,
first with independent heights and then given by a correlated harmonic
crystal; the upper tail of the independent wall decides whether the disorder
changes the leading repulsion height. When the disorder is in the field
itself, Chiarini and Pasqui \cite{ChiariniPasqui2025} obtained quenched
hard-wall asymptotics for uniformly elliptic random conductances under mixing
assumptions. Their rate involves the homogenized capacity and the essential
supremum of the on-site variances, and they also describe the conditional
profile and the local recentering.

Feldheim, Feldheim and Muirhead \cite{FeldheimFeldheimMuirhead2026} give a
broad account of persistence and entropic repulsion for stationary Gaussian
fields with a spectral singularity at the origin. Under regularity
assumptions, they express the logarithmic persistence cost and the repulsion
profile through capacity and equilibrium potentials. Their Section~1.6.3
also discusses moving walls; for the GFF, after normalization by
$\sqrt{\log N}$, it covers all levels strictly above the negative of the
leading repulsion height. The window studied here sits at the excluded
endpoint of this first-order regime, where the rate at speed
$N^{d-2}\log N$ vanishes.

\subsection{Setting and the critical window}\label{sec:setting}
Throughout the article, $d\ge3$. Let $g$ be the Green function of discrete-time simple random walk on
$\Z^d$, and let $\varphi=(\varphi_x)_{x\in\Z^d}$ be the centered Gaussian
field with covariance $g$, defined under a probability measure $\GP$. Put
$g_0=g(0,0)$. Let $D\subset\mathbb R^d$ be the closure of a nonempty
bounded Lipschitz domain and
set $D_N=(ND)\cap\Z^d$. 
We write $D^1(\R^d)$ for the homogeneous Sobolev
space with norm $\|\nabla f\|_2$ (see Section~\ref{sec:normalizations}).
As in \cite{LiSznitman2014,LiShiXu2024}, we use the Dirichlet form and the
capacity of Brownian motion:
\begin{equation}\label{eq:continuum-normalization}
\cE_{\R^d}(f,f)=\frac12\int_{\R^d}|\nabla f|^2\dd x,\qquad
 \capR(D)=\inf\left\{\cE_{\R^d}(f,f):
 f\in D^1(\R^d),\ f\ge1\text{ a.e. on }D\right\}.
\end{equation}
The remaining conventions are collected in
Section~\ref{sec:preliminaries}.

For $\beta\in\R$ we consider the level and the hard-wall event
\begin{align}
 h_N(\beta)&=\sqrt{g_0}\sqrt{4\log N+\log\log N+2\beta},
 \label{eq:window}\\
 \Omega_N(\beta)&=\{\varphi_x\ge-h_N(\beta)\text{ for every }x\in D_N\}.
 \label{eq:hardwall-event}
\end{align}
For the sequences $\beta_N$ considered below, $h_N(\beta_N)$ is well defined
and positive for large $N$. Throughout, we write
\begin{equation}\label{eq:aN}
 a_N=\frac{N^{d-2}}{\log N}
\end{equation}
for the speed of the critical window.
\paragraph{Heuristics.}
Why this window? Two costs have to be compared. 
\begin{itemize}
\item \emph{Local violations.} Write
$\theta=1/(2\sqrt{2\pi})$. By Mills' ratio (see \eqref{eq:mills}), a single site violates the
wall with probability
\[
 \GP\bigl(\varphi_0<-h_N(\beta)\bigr)\sim\frac{\theta e^{-\beta}}{N^2\log N},
\]
so the expected number of violations in $D_N$ is about
$\theta e^{-\beta}|D|\,a_N$. The term $\log\log N$ in
\eqref{eq:window} compensates the Gaussian prefactor and makes this number
exactly of order $a_N$. If the violations were independent, the
probability that none of them occurs would be about the exponential of minus
this number.
\item \emph{Collective displacement.} Lift the field by
$(g_0/h_N)f(\cdot/N)$, for a smooth bounded profile $f$. For bounded $s$, a
Taylor expansion of \eqref{eq:window} gives
\[
 h_N(\beta+s)-h_N(\beta)=\frac{g_0\,s}{h_N(\beta)}+O\bigl(h_N(\beta)^{-3}\bigr),
\]
so near $x$ the lift has the same effect as replacing $\beta$ by
$\beta+f(x/N)$ in \eqref{eq:window}. By the first item, it therefore
multiplies the violation probability at $x$ by $e^{-f(x/N)}(1+o(1))$. By the
Cameron--Martin formula and the energy scaling \eqref{eq:energy-scaling},
the lift costs
\[
 \frac{g_0^2}{2h_N^2}\,\cE_{\Z^d}\bigl(f(\cdot/N),f(\cdot/N)\bigr)
 \sim a_N\,\frac{g_0}{8d}\,\cE_{\R^d}(f,f).
\]
\end{itemize}
Both costs are of order $a_N$. Adding them suggests the
variational problem
\begin{equation}\label{eq:gff-variational}
 \cC_D^G(\beta)=\inf_{f\in D^1(\R^d)}
 \left\{\frac{g_0}{8d}\,\cE_{\R^d}(f,f)
       +\frac{e^{-\beta}}{2\sqrt{2\pi}}\int_D e^{-f(x)}\dd x\right\},
 \qquad \beta\in\R.
\end{equation}
The infimum may be restricted to $f\ge0$, since replacing $f$ by its
positive part lowers both terms. Our first result confirms the heuristic.

\begin{theorem}[Critical window]\label{thm:gff-window}
If $\beta_N\to\beta\in\R$, then
\begin{equation}\label{eq:critical-main}
 \lim_{N\to\infty}\frac1{a_N}
 \log\GP(\Omega_N(\beta_N))=-\cC_D^G(\beta).
\end{equation}
Equivalently, $a_N^{-1}\log\GP(\Omega_N(\beta))$ converges
to $-\cC_D^G(\beta)$ locally uniformly in $\beta\in\R$.
\end{theorem}

 We now list some basic properties of the variational problem.

\begin{proposition}[The optimal profile]\label{prop:soft-minimizer}
The variational problem \eqref{eq:gff-variational} has a unique minimizer
$f_\beta^G\in D^1(\R^d)$, and $f_\beta^G\ge0$. It solves, in the weak sense, the following differential equation
\begin{equation}\label{eq:soft-gff-euler}
 -\Delta f_\beta^G
 =\frac{4d}{\sqrt{2\pi}\,g_0}e^{-\beta-f_\beta^G}\mathbf 1_D.
\end{equation}
It has a bounded continuous representative, which is harmonic on
$\R^d\setminus\overline D$ and satisfies $f_\beta^G(x)=O(|x|^{2-d})$ as
$|x|\to\infty$.
\end{proposition}

\subsection{Outside the window}

The window separates two mechanisms. When $\beta_N\to-\infty$, the expected
number of violations is much larger than $a_N$. Suppressing them
forces a displacement of growing size, and the leading cost is given by the
capacity of $D$.

\begin{theorem}[Subcritical wall]\label{thm:subcritical-sharp}
If $\beta_N\to-\infty$ and $\beta_N^2=o(\log N)$, then
\begin{equation}\label{eq:subcritical-main}
 \lim_{N\to\infty}\frac{\log\GP(\Omega_N(\beta_N))}{a_N\beta_N^2}
 =-\frac{g_0}{8d}\capR(D).
\end{equation}
\end{theorem}

When $\beta_N\to+\infty$, the expected number of violations is
$o(a_N)$, so that a macroscopic displacement costs more than it
can save. The leading cost is then the expected number of violations.

\begin{theorem}[Supercritical wall]\label{thm:supercritical-sharp}
If $\beta_N\to+\infty$ and $\beta_N=o(\log N)$, then
\begin{equation}\label{eq:supercritical-main}
 \lim_{N\to\infty}\frac{\log\GP(\Omega_N(\beta_N))}{a_Ne^{-\beta_N}}
 =-\frac{|D|}{2\sqrt{2\pi}}.
\end{equation}
\end{theorem}

\begin{remark}[Matching of the three regimes]\label{rem:matching}
The off-window limits are those suggested by Theorem~\ref{thm:gff-window}:
\begin{equation}\label{eq:matching}
 \lim_{\beta\to-\infty}\frac{\cC_D^G(\beta)}{\beta^2}
 =\frac{g_0}{8d}\capR(D),\qquad
 \lim_{\beta\to+\infty}e^{\beta}\cC_D^G(\beta)=\frac{|D|}{2\sqrt{2\pi}}.
\end{equation}
For the second limit, write $\alpha=g_0/(8d)$ and $\lambda=\theta e^{-\beta}$.
The choice $f=0$ gives $\cC_D^G(\beta)\le\lambda|D|$. For $f\ge0$, the
bounds $e^{-f}\ge1-f$ and $\int_Df\le C_D\|\nabla f\|_2$ give
$\cC_D^G(\beta)\ge\lambda|D|-\lambda^2C_D^2/(2\alpha)$. For the first
limit, the upper bound follows by testing $f=tu$, with $u$ smooth, $u\ge1$
near $\overline D$, and $t=|\beta|+3\log|\beta|$. For the lower bound, let
$f_\beta$ be near-minimizers and $t=|\beta|$. The bound
$\cC_D^G(\beta)=O(t^2)$ forces
$|\{x\in D:f_\beta\le(1-\eps)t\}|=O(t^2e^{-\eps t})$, and the compactness
argument following \eqref{off:sub-upper} shows that every limit point of
$f_\beta/t$ is at least one almost everywhere on $D$. 
\end{remark}

\paragraph{Strategy of the proofs.}
By the domain Markov property, we split the field into independent
zero-boundary fields on mesoscopic boxes and a harmonic field determined by
the values on a finite skeleton (Section~\ref{sec:coarse}). The boxes have
volume of order $N^2\log N$, so that in the critical window each contains a
bounded expected number of violations. After a thin layer near the box
boundaries is discarded, a normal comparison inequality shows that the
violations inside one box are approximately Poisson, uniformly over subsets
of the box and over bounded, spatially varying shifts
(Section~\ref{sec:local}). This uniformity is what allows us to condition on
the harmonic field.

The rescaled harmonic field satisfies a large-deviation principle in local
$L^2$ with a Dirichlet-energy rate; the proof combines discrete energy
compactness, harmonic recovery sequences, and exponentially good
finite-dimensional approximations. In the critical window, the conditional
costs of the boxes add up to the exponential potential in
\eqref{eq:gff-variational}. The upper bound then follows from a Laplace
principle for a truncated potential, followed by removal of the truncation.
For the lower bound, a Cameron--Martin estimate forces the harmonic field
into a tube around a smooth profile, and Gaussian association restores the
discarded boundary layers (Section~\ref{sec:critical-proof}).

Outside the window the same scheme runs at moving speeds
(Section~\ref{sec:offwindow}). Below the window, if the rescaled harmonic
field lies below $(1-\eps)|\beta_N|$ on a region of macroscopic volume, then
the violations in that region cannot be suppressed at a cost of order
$a_N\beta_N^2$ (Lemma~\ref{off:low-regions}). The harmonic field must
therefore lie above a hard obstacle, and minimizing the
energy under this obstacle gives the capacity. Above the window, deviations
of the harmonic field live on a larger scale than the remaining violation
intensity, so the zero profile determines the leading Poisson cost.

\subsection{Interlacement coverage}

Random interlacements are the trace of a Poisson cloud of doubly infinite
transient nearest-neighbour trajectories on $\Z^d$, with an intensity
parameter $u\ge0$ \cite{Sznitman2007}. They describe in infinite volume the
local picture left by simple random walk on a large discrete torus in
dimension $d\ge3$ \cite{Sznitman2007,LiShiXu2024}. The cover level of a
finite set is the least intensity at which every point of the set has been
visited. For a box of side $N$ it is typically close to $dg_0\log N$, and
its fluctuations around this value were identified in \cite{Belius2012}. A
lower deviation asks for the whole box to be covered substantially earlier.
As for the GFF, a macroscopic increase of the occupation competes with the
suppression of the remaining vacant vertices.

Recall the definition of $D$ and $D_N$. 
Let $\mathcal I^u$ be the
interlacement set at level $u$, and write
$\mathcal T_N=\inf\{u:D_N\subset\mathcal I^u\}$ for the cover level of
$D_N$. We use the normalization $\Pp(x\notin\mathcal I^u)=e^{-u/g_0}$. At
level $\gamma dg_0\log N$, the expected number of vacant vertices in $D_N$
is of order $N^{d(1-\gamma)}$, which equals the capacitary scale $N^{d-2}$ when
$\gamma=2/d$. Below this value the cost of coverage is capacitary; at
$\gamma=2/d$ a smaller logarithmic speed appears.

\begin{proposition}\label{prop:ri-phase}
For every fixed $\gamma\in(0,2/d)$,
\begin{equation}\label{eq:ri-subcritical}
 \lim_{N\to\infty}\frac{\log\mathbb P(\mathcal T_N\leq\gamma dg_0\log N)}
 {N^{d-2}\log N}
 =-{g_0}\left(\sqrt\gamma-\sqrt{2/d}\right)^2
       \operatorname{cap}_{\mathbb R^d}(D).
\end{equation}
At $\gamma=2/d$,
\begin{align}
 -{\frac{g_0}{8d}}\operatorname{cap}_{\mathbb R^d}(D)
 &\leq\liminf_{N\to\infty}
       \frac{\log\mathbb P(\mathcal T_N\leq2g_0\log N)}{{a_N(\log\log N)^2}},
       \label{eq:ri-critical-lower}\\
 \limsup_{N\to\infty}
       \frac{\log\mathbb P(\mathcal T_N\leq2g_0\log N)}{{a_N(\log\log N)^2}}
 &\leq-{\frac{g_0}{32d}}\operatorname{cap}_{\mathbb R^d}(D).
       \label{eq:ri-critical-upper}
\end{align}
In particular, the speed at the transition is
${a_N(\log\log N)^2}$.
\end{proposition}

The change of regime at $\gamma=2/d$ was already described in
\cite[Section~0.1]{LiShiXu2024}: display~(0.10) there announces the
capacitary lower bound for $\gamma<2/d$, with the same constant as in
\eqref{eq:ri-subcritical}, and the paragraph after it
announces a matching upper bound to principal exponential order. Here we
give complete proofs of both bounds in \eqref{eq:ri-subcritical}, together
with the bounds at the transition. For comparison, the complementary regime
is described by \cite[(0.9)]{LiShiXu2024}:
\begin{equation}\label{eq:ri-supercritical}
 \lim_{N\to\infty}\frac{\log\mathbb P(\mathcal T_N\leq\gamma dg_0\log N)}
 {N^{d(1-\gamma)}}=-1,\qquad 2/d<\gamma<1.
\end{equation}
(Note that in \cite{LiShiXu2024} the statement is written for cubes only, however, there is no obstacle in extending the proof to general bounded Lipschitz domains.) Its lower bound follows from Poisson association, and the upper bound is
sketched in Appendix~B of that article. 

The proofs are given in Section~\ref{sec:ri}. For the upper bounds we
compare coverage with a GFF hard-wall event through Lupu's coupling
\cite{Lupu2016}. For the lower bounds, the entropy and vacancy identities of
\cite{LiSznitman2014} produce, at a controlled entropy cost, a configuration
with few vacant vertices, which an independent sprinkling then covers; this
avoids repeating the construction behind those identities. The constants in
\eqref{eq:ri-critical-lower} and \eqref{eq:ri-critical-upper} differ because
the RI and GFF windows have different logarithmic centerings.
Remark~\ref{rem:ri-window} explains this mismatch and predicts both the RI
window rate and the sharp constant at the transition.

\subsection{Discussion and open problems}\label{sec:discussion}

We list some questions raised by the results above.
\begin{enumerate}[label=(\arabic*)]
\item \emph{Entropic repulsion in the window.} The heuristic behind
\eqref{eq:gff-variational} suggests that, conditionally on
$\Omega_N(\beta)$, the field is lifted by about
$(g_0/h_N)f_\beta^G(\cdot/N)$, with $f_\beta^G$ the optimal profile of
Proposition~\ref{prop:soft-minimizer}. This lift is of order
$(\log N)^{-1/2}$, so, unlike the repulsion of
\cite{BDZ1995,DeuschelGiacomin1999}, it should not be visible in the local
law of the field; it should only be visible in the harmonic field $F_N$ of
Section~\ref{sec:coarse}. Does the interpolation $\widehat F_N$ converge in
probability to $f_\beta^G$ in $L^2_{\rm loc}(\R^d)$ under
$\GP(\,\cdot\mid\Omega_N(\beta))$? Uniqueness of the minimizer suggests a
route through the upper bound of Section~\ref{sec:critical-proof}, but the
removal of the truncation there would have to be made quantitative. In the
supercritical regime the analogous statement, with limit $0$, follows from
Proposition~\ref{prop:coarse-ldp} and \eqref{off:super-lower}: the harmonic
field deviates from $0$ only at speed $a_N$, while
$\log\GP(\Omega_N(\beta_N))$ is of order $a_Ne^{-\beta_N}=o(a_N)$. In the
subcritical regime one expects $\widehat F_N/|\beta_N|$ to converge to the
equilibrium potential of $D$.
\item \emph{The planar case.} For the two-dimensional Dirichlet field, the
hard-wall probability is described by the extremal process near the typical
minimum \cite{BiskupLouidor2020}, and by a relative capacity far from it
\cite{BolthausenDeuschelGiacomin2001,FelsLouidorWu2026}. The
double-exponential tail obtained in \cite{FelsLouidorWu2026} for the
conditioned minimum suggests that, for walls at height $s_N$ above the
typical minimum with $1\ll s_N\ll\log N$, the cost is purely capacitary, so
that no analogue of the window of Theorem~\ref{thm:gff-window} exists in
two dimensions. Confirming this, and identifying the constant, is open.
\item \emph{Sharp rates for interlacements.} Remark~\ref{rem:ri-window}
predicts the rate in the RI window \eqref{eq:ri-predicted-window} and the
sharp constant \eqref{eq:ri-critical-prediction} at the transition. A proof
would need an analogue of Lemma~\ref{lem:local-poisson} for vacant vertices,
uniform over slowly varying tilts of the local intensity. The normal
comparison inequality has no direct counterpart for interlacements, so the
correlations of the vacant set at mesoscopic distances would have to be
controlled by other means. A related question is the corresponding window
for the cover time of the discrete torus, in the setting of
\cite{LiShiXu2024}.

\item \emph{Random environment.} For the random-conductance model of
\cite{ChiariniPasqui2025}, the on-site variances fluctuate. Near the window,
where $h_N^2\approx4v\log N$ for a reference variance $v$, a site whose
variance is $v+\delta$ violates the wall with a probability larger by a
factor $N^{2\delta/v+O(\delta^2)}$. The window should
therefore be governed by the sites of nearly maximal variance, consistent
with the role of the essential supremum of the variances in
\cite{ChiariniPasqui2025}. Identifying the quenched window and its rate is
open.
\item \emph{Other Gaussian fields.} For the stationary fields of
\cite{FeldheimFeldheimMuirhead2026} with spectral density of order
$|\lambda|^{-\alpha}$ near the origin, $0<\alpha<d$, the GFF being the case
$\alpha=2$, the heuristic of Section~\ref{sec:setting} places the window at
$h_N^2=g_0(2\alpha\log N+\log\log N+2\beta)$, at speed
$N^{d-\alpha}/\log N$. The rate should then combine the exponential
potential with the homogeneous $\dot H^{\alpha/2}$ energy of the profile
in place of the Dirichlet energy.
\end{enumerate}

\subsection{Organization}
Section~\ref{sec:preliminaries} collects normalizations, mesoscopic scales
and Gaussian estimates, and proves Proposition~\ref{prop:soft-minimizer}.
Section~\ref{sec:local} proves the uniform local estimate.
Section~\ref{sec:coarse} establishes the large-deviation principle for the
harmonic field and its tube and Laplace consequences.
Section~\ref{sec:critical-proof} proves Theorem~\ref{thm:gff-window}, and
Section~\ref{sec:offwindow} treats the two off-window regimes.
Section~\ref{sec:ri} proves the interlacement results and compares the
logarithmic centerings of the RI and GFF windows.

\medskip

\noindent {\bf AI Usage Disclosure}: Commercial large language models were used to assist with the development of ideas (at certain stages of this project) and language polishing of the draft. The authors take full responsibility for the draft.

\medskip

\noindent {\bf Acknowledgements} The authors thank Pierre-Fran\c{c}ois Rodriguez for inspiring and encouraging discussions. X.L.~is supported by Beijing Natural Science Foundation (JQ26001) and the National Key R\&D Program of China (No.~2021YFA1002700).
J.S.\ is supported by an Overseas Bursary from Trinity College and a stipend from DPMMS.

\section{Preliminaries}\label{sec:preliminaries}

\subsection{Notation and normalizations}\label{sec:normalizations}

Positive constants $c,C$ may change from line to line. They may depend on
$d$ and $D$ without mention; any further dependence is indicated. We
recall $a_N=N^{d-2}/\log N$ from \eqref{eq:aN} and set
\begin{equation}\label{eq:basic-scales}
 \alpha=\frac{g_0}{8d},\qquad\theta=\frac1{2\sqrt{2\pi}}.
\end{equation}

\paragraph{Continuum.}
The homogeneous Sobolev space is
\[
 D^1(\R^d)=\{f\in L^{2d/(d-2)}(\R^d):\nabla f\in L^2(\R^d)\},
\]
normed by $\|\nabla f\|_2$; the Dirichlet form $\cE_{\R^d}$ and the
capacity $\capR$ are defined in \eqref{eq:continuum-normalization}. In the proofs we need
smooth competitors for the capacity \eqref{eq:continuum-normalization}. For a
Lipschitz domain they are available, and they show that
\eqref{eq:continuum-normalization} is also the capacity of the compact set
$\overline D$.

\begin{lemma}[Smooth competitors]\label{lem:smooth-capacity}
For a bounded Lipschitz domain $D$,
\begin{equation}\label{eq:smooth-capacity}
 \capR(D)=\inf\Bigl\{\cE_{\R^d}(h,h):\
 h\in C_c^\infty(\R^d),\ 0\le h\le1,\ h=1\text{ near }\overline D\Bigr\},
\end{equation}
where ``near $\overline D$'' means on some neighbourhood of $\overline D$.
\end{lemma}
This lemma follows from two standard facts about Newtonian
capacity. First, for a bounded Lipschitz domain, the interior
cone condition implies that the capacities of its interior
and closure coincide; see \cite{PortStone} for the
relevant potential-theoretic regularity results. Second,
Newtonian capacity admits a variational characterization
in terms of smooth compactly supported functions equal
to one on a neighbourhood of the compact set; see
\cite[Section 0]{NitzschnerSznitman2020}.
Together with the Sobolev variational characterization
of capacity, these facts yield \eqref{eq:smooth-capacity}.

\paragraph{Lattice.}
Let $P$ be the transition operator of discrete-time simple random walk, so
that $g=(I-P)^{-1}$. The corresponding lattice energy is
\begin{equation}\label{eq:lattice-energy}
 \cE_{\Z^d}(u,u)=\langle u,(I-P)u\rangle
 =\frac1{4d}\sum_{x\in\Z^d}\sum_{|e|=1}(u(x+e)-u(x))^2.
\end{equation}
The density of $\varphi$ is formally proportional to
$\exp\{-\frac12\cE_{\Z^d}(\varphi,\varphi)\}$, so a deterministic shift
$m$ has Cameron--Martin cost $\frac12\cE_{\Z^d}(m,m)$. For
$f\in C_c^\infty(\R^d)$, a Taylor expansion gives
\begin{equation}\label{eq:energy-scaling}
 N^{2-d}\cE_{\Z^d}(f(\cdot/N),f(\cdot/N))
 \longrightarrow\frac1{2d}\int_{\R^d}|\nabla f|^2\dd x=\frac1d\cE_{\R^d}(f,f).
\end{equation}
The lattice capacity of a finite set $K$ is
$\capZ(K)=\inf\{\cE_{\Z^d}(u,u):u\ge1\text{ on }K\}$, the infimum being
over finitely supported $u$. Its scaling limit is classical: for a bounded
Lipschitz domain $D$,
\begin{equation}\label{eq:capacity-scaling}
 N^{2-d}\capZ(D_N)\longrightarrow\frac1d\capR(D);
\end{equation}
see \cite[Proposition~2.3]{Bouchot2025} and \cite{LiSznitman2014}. In
\cite{Bouchot2025} the limit is stated for $(NK)\cap\Z^d$ with $K$ compact
and regular, that is, $\capR(K)=\capR(\operatorname{Int}K)$. It applies to
$K=\overline D$, which is regular by Lemma~\ref{lem:smooth-capacity}, and
since $(NK)\cap\Z^d\subset D_N\subset(N\overline D)\cap\Z^d$ for compact
regular $K\subset D$ (for instance finite unions of closed grid cubes),
\eqref{eq:capacity-scaling} follows by monotonicity and inner regularity of
capacity.

\subsection{Boxes, scales and Green bounds}\label{sec:scales}

We fix once and for all
\begin{equation}\label{eq:mesoscopic-scales}
 \ell_N=\lfloor(N^2\log N)^{1/d}\rfloor,\qquad
 r_N=(\log N)^{4/(d-2)}.
\end{equation}
The side $\ell_N$ is chosen so that $\ell_N^d\sim N^2\log N$: in the
critical window, a box of side $\ell_N$ then contains an expected number
of violations that tends to $\theta e^{-\beta}$; in particular, it is bounded
and bounded away from zero. We also have $\ell_N/N\to0$, which is used for
harmonic recovery. The trimming distance $r_N$ satisfies
\begin{equation}\label{eq:trimming-conditions}
 (\log N)^2r_N^{2-d}=(\log N)^{-2}\longrightarrow0,
 \qquad (\log N)\frac{r_N}{\ell_N}=N^{-2/d+o(1)}\longrightarrow0.
\end{equation}
The first condition makes the variance of the harmonic field negligible at
distance $r_N$ from the box boundaries, even after multiplication by $\log
N$. The second makes the discarded boundary layers small. The polylogarithmic
size of $r_N$ is used once more in Section~\ref{sec:offwindow}, to control
these layers at moving speeds.

A \emph{box} is a set $B=z+\{1,\dots,\ell_N-1\}^d$ with $z\in\Z^d$. Its
zero-boundary GFF $\eta^B$ is the centered Gaussian field on $B$ with
covariance $g_B$, the Green function of the walk killed on leaving $B$, and
its trimmed version is
\[
 B^-=\{x\in B:d_\infty(x,B^c)\ge r_N\}.
\]
The standard Green bounds \cite[Theorem~4.3.1]{LawlerLimic2010} give
\begin{equation}\label{eq:green-bounds}
 0\le g(x,y)\le C(1+|x-y|)^{2-d},\qquad
 0\le g_0-g_B(x,x)\le C r^{2-d}
 \quad\text{if }d_\infty(x,B^c)\ge r.
\end{equation}

In particular, the variances of the field $\eta^B$ on $B^-$ are
$g_0+O((\log N)^{-4})$. We also use that
\begin{equation}\label{eq:green-offdiagonal}
 \sup_{x\ne y}g(x,y)=g(0,e_1)=g_0-1.
\end{equation}

\subsection{Gaussian tools}\label{sec:gaussian-facts}

We write $\Phi$ and $\phi$ for the distribution function and the density of
a standard normal variable. Mills' ratio gives
\begin{equation}\label{eq:mills}
 1-\Phi(t)=(1+o(1))\frac{\phi(t)}t\qquad(t\to\infty);
\end{equation}
see \cite[(7.8.2)]{DLMF}. 
Gaussian vectors with nonnegative covariances are associated \cite{Pitt1982}: increasing events are positively correlated. We
use the following form of the normal comparison inequality of Li and Shao
\cite[Corollary~2.1]{LiShao2002}. If $(X_i)$ is a centered Gaussian vector
with unit variances and correlations $\rho_{ij}$, then
\begin{equation}\label{eq:normal-comparison}
 \left|\Pp(X_i\le u_i,\ \forall i)-\prod_i\Phi(u_i)\right|
 \le\frac14\sum_{i<j}|\rho_{ij}|
       \exp\left\{-\frac{u_i^2+u_j^2}{2(1+|\rho_{ij}|)}\right\}.
\end{equation}

\subsection{The optimal profile}\label{sec:variational}\label{sec:soft-obstacle}

We prove Proposition~\ref{prop:soft-minimizer} for the more general
functional
\begin{equation}\label{eq:soft-functional}
 J_{a,\lambda}(f)
 =a\,\cE_{\R^d}(f,f)
   +\lambda\int_D e^{-f(x)}\dd x,
 \qquad f\in D^1(\R^d),\quad a,\lambda>0,
\end{equation}
with values in $[0,\infty]$. The proposition is the case $a=\alpha$ and
$\lambda=\theta e^{-\beta}$.

\begin{proof}[Proof of Proposition~\ref{prop:soft-minimizer}]
\emph{Existence.} Since $J_{a,\lambda}(f_+)\le J_{a,\lambda}(f)$ and
$J_{a,\lambda}(0)=\lambda|D|$, there is a nonnegative minimizing sequence,
and it is bounded in $D^1$. By the Sobolev inequality and local Rellich
compactness, a subsequence converges weakly in $D^1$ and almost everywhere.
Weak lower semicontinuity of the energy (see, e.g.,
\cite[Proposition~3.5(iii), p.~58]{Brezis2011}) and dominated convergence on $D$
(note $e^{-f}\le1$ for $f\ge0$) show that the limit is a minimizer.

\emph{Uniqueness and sign.} The energy is strictly convex on $D^1$ and the
potential term is convex, so the minimizer is unique; call it
$f_{a,\lambda}$. If its negative part $f_-$ did not vanish, then
$\|\nabla f_-\|_2>0$ and
$J_{a,\lambda}(f_+)\le J_{a,\lambda}(f_{a,\lambda})-\tfrac a2\|\nabla f_-\|_2^2$,
which is impossible. Hence $f_{a,\lambda}\ge0$.

\smallskip
\noindent \emph{Euler--Lagrange equation.} For $\psi\in C_c^\infty$ and $|t|\le1$ we
have $e^{-f_{a,\lambda}-t\psi}\le e^{\|\psi\|_\infty}$, so we may
differentiate under the integral at $t=0$:
\begin{equation}\label{eq:soft-euler}
 a\int_{\R^d}\nabla f_{a,\lambda}\cdot\nabla\psi\dd x
 =\lambda\int_D e^{-f_{a,\lambda}}\psi\dd x.
\end{equation}
This extends to $\psi\in D^1$ by density and the bound
$\|\psi\|_{L^1(D)}\le C_D\|\nabla\psi\|_2$, which follows from H\"older's
inequality on $D$ and the Sobolev inequality
$\|\psi\|_{L^{2d/(d-2)}}\le C_d\|\nabla\psi\|_2$. With $a=\alpha$ and
$\lambda=\theta e^{-\beta}$ we have $\lambda/a=4de^{-\beta}/(\sqrt{2\pi}\,g_0)$,
which is \eqref{eq:soft-gff-euler}.

\smallskip
\noindent \emph{Regularity.} Let $\Gamma(x)=((d-2)|\mathbb S^{d-1}|)^{-1}|x|^{2-d}$ be
the fundamental solution of $-\Delta$. The source
$\lambda e^{-f_{a,\lambda}}\mathbf1_D/{a}$ is bounded and compactly
supported, so its Newtonian potential belongs to $D^1$ and has the same weak
Laplacian as $f_{a,\lambda}$. Testing the difference of the two functions
against itself shows that they coincide:
\begin{equation}\label{eq:soft-potential}
 f_{a,\lambda}(x)
 ={\frac{\lambda}{a}}\int_D\Gamma(x-y)e^{-f_{a,\lambda}(y)}\dd y
 \le{\frac{\lambda}{a}}\int_D\Gamma(x-y)\dd y.
\end{equation}
This representation gives continuity, boundedness, harmonicity off
$\overline D$, and $f_{a,\lambda}(x)=O(|x|^{2-d})$ at infinity.
\end{proof}

\section{A uniform local extreme-value estimate}\label{sec:local}

Recall conventions regarding boxes and scales in Section \ref{sec:scales}. For a given box $B$ and a zero-boundary GFF $\eta^B$ inside $B$, we show that, away
from the box boundary, the probability that it stays above a slightly shifted
wall is the one predicted by independent sites, with a relative error that is
uniform over all subsets of the box and all bounded shifts. This estimate is
the input common to the three regimes of
Theorems~\ref{thm:gff-window}--\ref{thm:supercritical-sharp}.

We start with a one-point bound.
\begin{lemma}\label{lem:one-point}
For $Z\sim N(0,g_0)$, uniformly for $\beta$ and $q$ in fixed compact sets,
\begin{equation}\label{eq:one-point}
 \Pp\left(Z<-h_N(\beta)-\frac{g_0q}{h_N(\beta)}\right)
 =(1+o(1))\frac{\theta e^{-\beta-q}}{N^2\log N}.
\end{equation}
The same holds for $Z\sim N(0,g_B(x,x))$, uniformly over boxes $B$ and
$x\in B^-$.
\end{lemma}
\begin{proof}
Put $u=h_N+g_0q/h_N$. By Mills' ratio \eqref{eq:mills},
\[
 \Pp(Z<-u)=(1+o(1))\frac{\sqrt{g_0}}{\sqrt{2\pi}\,u}
 e^{-u^2/(2g_0)},
 \qquad
 \frac{u^2}{2g_0}=2\log N+\tfrac12\log\log N+\beta+q
                 +\frac{g_0q^2}{2h_N^2}.
\]
The last term is $o(1)$, and the prefactor is
$(1+o(1))/(2\sqrt{2\pi\log N})$; this proves \eqref{eq:one-point}. 
Replacing
$g_0$ by $g_B(x,x)$ changes the exponent by $O(h_N^2r_N^{2-d})=o(1)$, by
\eqref{eq:green-bounds}.
\end{proof}

\begin{lemma}[Comparison at the mesoscopic scale]\label{lem:mesoscopic-comparison}
There is $c=c(d)>0$ with the following property. Uniformly over boxes $B$,
subsets $A\subset B^-$, and positive thresholds $(z_x)_{x\in A}$ with
$z_x^2=4\log N+o(\log N)$ uniformly in $x\in A$,
\begin{equation}\label{eq:mesoscopic-comparison}
 \Delta_A:=\left|\Pp\left(-\frac{\eta_x^B}{\sqrt{g_B(x,x)}}\le z_x,
                  \ x\in A\right)-\prod_{x\in A}\Phi(z_x)\right|
 \le N^{-c+o(1)}.
\end{equation}
If moreover $z_x^2=4\log N+\log\log N+O(1)$ uniformly, then
\begin{equation}\label{eq:comparison-relative}
 \Delta_A=o\left(\frac{|A|}{N^2\log N}\right).
\end{equation}
All error terms are uniform as long as the assumptions on the thresholds hold
uniformly; for $A=\varnothing$, $\Delta_A=0$.
\end{lemma}
\begin{proof}
Let $\rho_{xy}$ denote the correlations of the standardized field
$\eta_x^B/\sqrt{g_B(x,x)}$, $x\in A$; the field $-\eta^B$ has the same
correlations. They are nonnegative, and since $g_B\le g$,
\eqref{eq:green-bounds} gives $\rho_{xy}\le C|x-y|^{2-d}$. Moreover
$g_B(x,x)\ge g_0-Cr_N^{2-d}$ on $B^-$, so by \eqref{eq:green-offdiagonal},
for all large $N$,
\[
 \rho_{xy}\le\rho_*:=1-\frac1{2g_0}<1\qquad(x\ne y\in A).
\]
We apply \eqref{eq:normal-comparison} and split the sum at the distance
$\varrho_N=(\log N)^\kappa$, where $\kappa>1/(d-2)$ is fixed.

\smallskip
\noindent \emph{Near pairs.} There are at most $C|A|\varrho_N^d$ pairs with
$|x-y|\le\varrho_N$, and each contributes at most
$\exp\{-(z_x^2+z_y^2)/(2(1+\rho_*))\}=N^{-4/(1+\rho_*)+o(1)}$.

\smallskip
\noindent \emph{Far pairs.} If $|x-y|>\varrho_N$, then
$\rho_{xy}\le C(\log N)^{-\kappa(d-2)}$, so $\rho_{xy}\log N\to0$ and the
exponential factor is $N^{-4+o(1)}$. Since
$\sum_{y\in B\setminus\{x\}}|x-y|^{2-d}\le C\ell_N^2$, the far pairs
contribute at most $C|A|\ell_N^2N^{-4+o(1)}$.

Since $|A|\le\ell_N^d=N^{2+o(1)}$ and $\ell_N^2=N^{4/d+o(1)}$, we obtain
\eqref{eq:mesoscopic-comparison} with any fixed
\[
 0<c<\min\left\{\frac4{1+\rho_*}-2,\ 2-\frac4d\right\}.
\]
Under the sharper assumption on the thresholds, the far-pair bound improves
to $C|A|\ell_N^2N^{-4}(\log N)^{-1}$. Dividing it, and the near-pair bound,
by $|A|/(N^2\log N)$ gives $O(\ell_N^2/N^2)$ and $N^{2-4/(1+\rho_*)+o(1)}$
respectively. Both tend to zero, which proves \eqref{eq:comparison-relative}.
\end{proof}

\begin{lemma}[Uniform local Poisson cost]\label{lem:local-poisson}
Fix $M<\infty$ and a compact interval $J\subset\R$. Uniformly over
$\beta\in J$, boxes $B$, sets $A\subset B^-$, and functions
$q:A\to[-M,M]$,
\begin{equation}\label{eq:local-poisson}
 \log\Pp\left(\eta_x^B\ge-h_N(\beta)-\frac{g_0}{h_N(\beta)}q_x,
                    \ x\in A\right)
 =-(1+o(1))\frac{\theta e^{-\beta}}{N^2\log N}
                    \sum_{x\in A}e^{-q_x}.
\end{equation}
\end{lemma}
\begin{proof}
Let $p_x$ be the probability that site $x$ violates its threshold, and
$\mu_A=\sum_{x\in A}p_x$. By Lemma~\ref{lem:one-point},
\begin{equation}\label{eq:local-intensity}
 \mu_A=(1+o(1))\frac{\theta e^{-\beta}}{N^2\log N}
                 \sum_{x\in A}e^{-q_x},\qquad
 \mu_A\asymp\frac{|A|}{N^2\log N}=O(1).
\end{equation}
The standardized thresholds satisfy $z_x^2=4\log N+\log\log N+O(1)$
uniformly in the parameters, so Lemma~\ref{lem:mesoscopic-comparison} bounds
the additive comparison error by $o(\mu_A)$, also when $\mu_A\to0$. For
independent sites,
\[
 \log\prod_{x\in A}(1-p_x)=-\mu_A+O\bigl(\mu_A\max_xp_x\bigr)
                         =-\mu_A+o(\mu_A).
\]
This product is bounded away from zero by \eqref{eq:local-intensity}, so the
additive error $o(\mu_A)$ also gives an error $o(\mu_A)$ after taking
logarithms.
\end{proof}

\section{The macroscopic harmonic field}\label{sec:coarse}

We now introduce the harmonic part of the field and prove a large-deviation
principle for it. The two ingredients are energy compactness and an
exponentially good approximation by finite-dimensional projections.
Throughout this section, $h_N>0$ is deterministic and $h_N^2\sim4g_0\log N$;
for instance, $h_N=h_N(\beta_N)$ with $\beta_N=o(\log N)$.

\subsection{The skeleton and the harmonic field}

Let $L_N=\ell_N\lceil N^2/\ell_N\rceil$, $K_N=[-L_N,L_N]^d\cap\Z^d$, and
define the skeleton
\begin{equation}\label{eq:skeleton}
 \Sigma_N=\{x\in K_N:\ x_i\in\ell_N\Z\text{ for some }i\}.
\end{equation}
Since $L_N$ is a multiple of $\ell_N$, the set $K_N\setminus\Sigma_N$ is the
disjoint union of the boxes $z+\{1,\dots,\ell_N-1\}^d$ with $z\in\ell_N\Z^d$
that are contained in $K_N$; these are the boxes of the decomposition. With
the conventions of Section~\ref{sec:scales}, a box $B$ is the open lattice
interior, its outer vertex boundary $\{y\notin B:\ |y-x|=1\text{ for some
}x\in B\}$ is contained in $\Sigma_N$, and $\eta^B$ is the field killed on
leaving $B$, that is, on $\Sigma_N$. Two distinct boxes contain no
neighbouring vertices, since a layer of $\Sigma_N$ separates them. Since
$L_N/N\to\infty$, the sets $K_N/N$ exhaust $\R^d$. 

Define, on the whole lattice,
\begin{equation}\label{eq:coarse-definition}
 \cG_N=\sigma(\varphi_z:z\in\Sigma_N),\qquad
 m_N(x)=\E[\varphi_x\mid\cG_N],\qquad
 F_N(x)=\frac{h_N}{g_0}m_N(x).
\end{equation}
By the domain Markov property (see \cite{Sheffield2007}),
$\varphi_x=m_N(x)+\eta_x^B$ for $x$ in a box $B$, where the fields $\eta^B$
are independent zero-boundary GFFs, independent of $\cG_N$. In this
decomposition, the hard-wall constraint at $x\in B$ reads
\[
 \varphi_x\ge-h_N\iff\eta_x^B\ge-h_N-\frac{g_0}{h_N}F_N(x),
\]
which is the form appearing in Lemma~\ref{lem:local-poisson}, with
$q_x=F_N(x)$.

Writing $G_{\Sigma_N}$ for the Green matrix on $\Sigma_N$, we have
$m_N(x)=g(x,\Sigma_N)G_{\Sigma_N}^{-1}\varphi_{\Sigma_N}$, which is the
unique extension of its skeleton values which is harmonic on
$\Z^d\setminus\Sigma_N$ and belongs to $\ell^{2^*}(\Z^d)$, where
$2^*=2d/(d-2)$. In particular, its values outside $K_N$ are determined in
this way. Let $\cH_N$ denote the finite-dimensional space of such
extensions. With the normalization \eqref{eq:lattice-energy}, the
Cameron--Martin cost of the law of $F_N$ is
\begin{equation}\label{eq:coarse-CM}
 I_N(u)=\frac{g_0^2}{2h_N^2}\cE_{\Z^d}(u,u),
 \qquad u\in\cH_N.
\end{equation}
Indeed, the precision operator of the full field is $I-P$, and minimizing
the energy over extensions of prescribed skeleton values (equivalently,
taking the Gaussian Schur complement) gives \eqref{eq:coarse-CM}.

\paragraph{Trimmed sets.}
Let $D_N^-$ be the set of sites of $D_N$ that lie in the trimmed part $B^-$
of some box $B$ of the decomposition, and put $R_N=D_N\setminus D_N^-$.
Since $\ell_N^d\sim N^2\log N$, at most $Ca_N$ boxes meet $D_N$, and
\begin{equation}\label{eq:RN-size}
 |R_N|\le CN^d\frac{r_N}{\ell_N}=o(N^d).
\end{equation}

\paragraph{Interpolation and rates.}
For a lattice function $u$, write $\widehat u_N$ for its rescaled piecewise-constant
interpolation, equal to $u(x)$ on $Q_{N,x}=x/N+[0,N^{-1})^d$. Recall \eqref{eq:basic-scales} and set
\begin{equation}\label{eq:coarse-q-rate}
 q_N(u)=\frac{I_N(u)}{a_N},\qquad
 I(f)=
 \begin{cases}
  \alpha\cE_{\R^d}(f,f),&f\in D^1(\R^d),\\
  +\infty,&\text{otherwise}.
 \end{cases}
\end{equation}

\subsection{Energy compactness and recovery}

The next lemma identifies $I$ as the limit of the rescaled
Cameron--Martin costs $q_N$.

\begin{lemma}\label{lem:discrete-energy}
Suppose $u_N\in\cH_N$ and $\sup_Nq_N(u_N)<\infty$. Then $(\widehat u_N)$ is
precompact in $L^2_{\mathrm{loc}}(\R^d)$, and if
$\widehat u_{N_j}\to f$ in $L^2_{\mathrm{loc}}$, then
\begin{equation}\label{eq:coarse-liminf}
 I(f)\le\liminf_j q_{N_j}(u_{N_j}).
\end{equation}
Conversely, for every $f\in C_c^\infty(\R^d)$ there are $u_N^f\in\cH_N$ with
\begin{equation}\label{eq:coarse-recovery}
 \sup_{x\in\Z^d}|u_N^f(x)-f(x/N)|\longrightarrow0,
 \qquad
 \widehat u_N^f\longrightarrow f\ \text{in }L^2_{\mathrm{loc}},
 \qquad q_N(u_N^f)\longrightarrow I(f).
\end{equation}
\end{lemma}

\begin{proof}
\emph{Interpolation bounds.} Let $\widetilde u_N$ be the multilinear
interpolation of $u$ {at scale $1/N$: on each closed cube
$\overline{Q_{N,x}}$, it is the multilinear function that equals $u(y)$ at
the corner $y/N$, for $y\in x+\{0,1\}^d$. We claim that, for every
$u\in\ell^{2^*}(\Z^d)$ of finite energy,}
\begin{align}
 \|\widehat u_N\|_{L^{2^*}(\R^d)}^2
   &\le C_dN^{2-d}\cE_{\Z^d}(u,u),\label{eq:coarse-sobolev}\\
 \|\nabla\widetilde u_N\|_{L^2(\R^d)}^2
   &\le2dN^{2-d}\cE_{\Z^d}(u,u),\label{eq:coarse-interpolation}\\
 \|\widetilde u_N-\widehat u_N\|_{L^2(\R^d)}^2
   &\le C_dN^{-d}\cE_{\Z^d}(u,u).\label{eq:coarse-interpolation-error}
\end{align}
{We first assume that $u$ is finitely supported.} The partial derivatives
of $\widetilde u_N$ are $N$ times convex
combinations of forward differences, so Jensen's inequality gives
\eqref{eq:coarse-interpolation}, and summing along edges gives
\eqref{eq:coarse-interpolation-error}. For \eqref{eq:coarse-sobolev},
fix a cube $Q_{N,x}$. After rescaling it to $[0,1]^d$, the quantity
$(N^d\int_{Q_{N,x}}|\widetilde u_N|^{2^*})^{1/2^*}$ is a norm of the $2^d$
corner values $(u(y))_{y\in x+\{0,1\}^d}$: a multilinear function that
vanishes a.e.\ on $[0,1]^d$ vanishes at the corners. By equivalence of norms
on $\R^{2^d}$, the corner value $u(x)$ is controlled by this norm, that is,
\[
 \int_{Q_{N,x}}|\widehat u_N|^{2^*}=N^{-d}|u(x)|^{2^*}
 \le C_d\int_{Q_{N,x}}|\widetilde u_N|^{2^*}.
\]
Summing over $x$, applying the Sobolev inequality on $\R^d$ to the compactly
supported Lipschitz function $\widetilde u_N$, and then using
\eqref{eq:coarse-interpolation}, we obtain
\[
 \|\widehat u_N\|_{L^{2^*}(\R^d)}^2\le C_d\|\widetilde u_N\|_{L^{2^*}(\R^d)}^2
 \le C_d\|\nabla\widetilde u_N\|_{L^2(\R^d)}^2\le C_dN^{2-d}\cE_{\Z^d}(u,u),
\]
which is \eqref{eq:coarse-sobolev}. For a general $u\in\ell^{2^*}(\Z^d)$ of
finite energy, let $u_n(x)=u(x)\chi(x/n)$ with $\chi\in C_c^\infty(\R^d)$
equal to one near the origin. Then $\cE_{\Z^d}(u-u_n,u-u_n)\to0$: the
gradient of $u-u_n$ is bounded by that of $u$ away from the origin plus
$|u|\,|\nabla\chi(\cdot/n)|/n$, and the latter is controlled by H\"older's
inequality and $u\in\ell^{2^*}$. Since $u_n\to u$ pointwise, Fatou's lemma
transfers the three bounds from $u_n$ to $u$.

\smallskip
\noindent \emph{Compactness.} By \eqref{eq:coarse-CM} and $h_N^2\sim4g_0\log N$,
\begin{equation}\label{eq:coarse-q-asymptotic}
 q_N(u)=\left(\frac{g_0}{8}+o(1)\right)
                 N^{2-d}\cE_{\Z^d}(u,u).
\end{equation}
A bound on $q_N$ thus bounds $\widetilde u_N$ in $H^1(B_R)$ for every ball
$B_R$, by \eqref{eq:coarse-sobolev}--\eqref{eq:coarse-interpolation-error}.
Rellich compactness \cite[Sections~5.6--5.7]{Evans2010}, together with
\eqref{eq:coarse-interpolation-error}, gives precompactness of
$\widehat u_N$ in $L^2_{\rm loc}$.

\smallskip
\noindent \emph{Lower bound.} To keep the sharp constant, we use difference quotients
rather than $\nabla\widetilde u_N$. Define $D_{N,i}u=N(u(x+e_i)-u(x))$ on
$Q_{N,x}$. Then
\begin{equation}\label{eq:coarse-difference-energy}
 \sum_{i=1}^d\|D_{N,i}u\|_{L^2(\R^d)}^2
       =2dN^{2-d}\cE_{\Z^d}(u,u).
\end{equation}
Summation by parts against smooth compactly supported test functions shows
that every weak limit of $D_{N_j,i}u_{N_j}$ equals $\partial_if$. By weak lower
semicontinuity and \eqref{eq:coarse-difference-energy},
\[
 {\int_{\R^d}|\nabla f|^2\dd x}
 \le2d\liminf_jN_j^{2-d}\cE_{\Z^d}(u_{N_j},u_{N_j}).
\]
The global $L^{2^*}$ bound passes to the limit, so $f\in D^1(\R^d)$.
Together with \eqref{eq:coarse-q-asymptotic} and $\alpha d=g_0/8$ (recall
$\cE_{\R^d}(f,f)=\frac12\int|\nabla f|^2$), this proves \eqref{eq:coarse-liminf}.

\smallskip
\noindent \emph{Recovery.} Prescribe the values $f(x/N)$ on $\Sigma_N$ and let $u_N^f$
be their canonical harmonic extension. Since $K_N/N$ exhausts $\R^d$, the
data on the outer boundary of $K_N$ vanish for large $N$, and then $u_N^f=0$
outside $K_N$. Since $f(\cdot/N)$ is another finite-energy extension of the
same data, energy minimization and \eqref{eq:energy-scaling} give
\begin{equation}\label{eq:coarse-recovery-energy}
 \cE_{\Z^d}(u_N^f,u_N^f)
 \le\cE_{\Z^d}(f(\cdot/N),f(\cdot/N)),\qquad
 N^{2-d}\cE_{\Z^d}(f(\cdot/N),f(\cdot/N))
 \longrightarrow\frac1d\cE_{\R^d}(f,f).
\end{equation}
In a box, the difference $u_N^f-f(\cdot/N)$ vanishes on the boundary and has
discrete Laplacian bounded by $C_fN^{-2}$. Its random-walk representation and
the bound on expected exit time
\begin{equation}\label{eq:coarse-recovery-uniform}
 \sup_x|u_N^f(x)-f(x/N)|\le C_f\frac{\ell_N^2}{N^2}=o(1).
\end{equation}
Only boxes within distance $O(\ell_N/N)$ of $\supp f$ (in macroscopic units)
carry a nonzero $u_N^f$, so $\widehat u_N^f\to f$ in $L^2_{\rm loc}$. The
energy bound \eqref{eq:coarse-recovery-energy} gives
$\limsup_Nq_N(u_N^f)\le I(f)$, and \eqref{eq:coarse-liminf} gives the
matching lower bound.
\end{proof}

\subsection{Small fluctuations and shifted tubes}

Away from the skeleton, the harmonic field has small variance. Let $x\in B^-$
for a box $B$ of the decomposition (for instance $x\in D_N^-$). By the
domain Markov decomposition and
\eqref{eq:green-bounds},
\[
 \Var(m_N(x))=g_0-g_B(x,x)\le C_dr_N^{2-d},
 \qquad\text{hence}\qquad
 \Var(F_N(x))\le C(\log N)r_N^{2-d}.
\]
At every site, $\Var(F_N(x))\le h_N^2/g_0\le C\log N$. Within a fixed ball
$B_R$, the rescaled sites at distance less than $r_N$ from the skeleton
occupy a volume at most $C_R(r_N/\ell_N+N^{-1})$. Consequently, by
\eqref{eq:trimming-conditions},
\begin{equation}\label{eq:coarse-vanishing-trace}
 \E\|\widehat F_N\|_{L^2(B_R)}^2
 \le C_R\log N\left(r_N^{2-d}+\frac{r_N}{\ell_N}+N^{-1}\right)
 \longrightarrow0.
\end{equation}
A Gaussian tail bound and a union bound over $D_N^-$ also give, for every
fixed $\eps>0$,
\begin{equation}\label{eq:coarse-centered-sup}
 \GP\left(\max_{x\in D_N^-}|F_N(x)|>\eps\right)
 \le2|D_N|\exp\left\{-\frac{c\eps^2}{(\log N)r_N^{2-d}}\right\}
 \longrightarrow0,
\end{equation}
since $(\log N)r_N^{2-d}=(\log N)^{-3}$ by \eqref{eq:mesoscopic-scales}, so that the exponent is of order
$\eps^2(\log N)^3$, while $\log(2|D_N|)=O(\log N)$.

\begin{lemma}\label{lem:tube}
For a symmetric Borel set $A\subset\cH_N$ and $u\in\cH_N$,
\begin{equation}\label{eq:coarse-CM-tube}
 \GP(F_N\in u+A)\ge e^{-I_N(u)}\GP(F_N\in A).
\end{equation}
Let $t_N\ge1$ be deterministic. For $f\in C_c^\infty(\R^d)$ and
$\eps>0$,
\begin{equation}\label{eq:coarse-shifted-tube}
 \liminf_{N\to\infty}\frac1{a_Nt_N^2}
 \log\GP\left(\max_{x\in D_N^-}
       |F_N(x)/t_N-f(x/N)|\le\eps\right)\ge-I(f).
\end{equation}
\end{lemma}

\begin{proof}
By the Cameron--Martin formula, the translated probability equals
$e^{-I_N(u)}$ times the integral over $A$ of the exponential of a linear
functional. By symmetry of $A$, this exponential may be replaced by its
hyperbolic cosine, which is at least one; see also
\cite[Chapter~2]{Bogachev1998}.

For \eqref{eq:coarse-shifted-tube}, apply \eqref{eq:coarse-CM-tube} with
$u=t_Nu_N^f$ and $A=\{v\in\cH_N:\max_{D_N^-}|v|\le\eps t_N/2\}$, where
$u_N^f$ is introduced in Lemma~\ref{lem:discrete-energy}. By
\eqref{eq:coarse-recovery-uniform}, the event $\{F_N\in t_Nu_N^f+A\}$ is
contained in the event in \eqref{eq:coarse-shifted-tube} for large $N$, and
by \eqref{eq:coarse-centered-sup} and $t_N\ge1$, $\GP(F_N\in A)\to1$. Since
$I_N(t_Nu_N^f)=t_N^2a_Nq_N(u_N^f)$ by \eqref{eq:coarse-CM} and
\eqref{eq:coarse-q-rate}, \eqref{eq:coarse-recovery} concludes the
proof.
\end{proof}

\subsection{The Gaussian large-deviation principle}

\begin{proposition}\label{prop:coarse-ldp}
The laws of $\widehat F_N$ satisfy an LDP on $L^2_{\mathrm{loc}}(\R^d)$, with
its usual topology, at speed $a_N$ and with good rate function $I$ from
\eqref{eq:coarse-q-rate}. More generally, for every deterministic sequence
$t_N\ge1$, the laws of $\widehat F_N/t_N$ satisfy the LDP at speed $a_Nt_N^2$
with the same rate function.
\end{proposition}

\begin{proof}
We prove the general statement; set $X_N=F_N/t_N$ and $a_N^{(t)}=a_Nt_N^2$.
The law of $X_N$ has Cameron--Martin cost $a_N^{(t)}q_N$.

\smallskip
\noindent \emph{Goodness.} A sublevel set of $I$ is bounded in $L^{2^*}(\R^d)$ and in
$H^1(B_R)$ for every $R$, so Rellich compactness and weak lower
semicontinuity show that it is compact in $L^2_{\rm loc}$.

\smallskip
\noindent \emph{Finite-dimensional projections.} For
$\psi_1,\ldots,\psi_k\in C_c^\infty(\R^d)$, let
$p_\psi(f)=(\int f\psi_j)_{j=1}^k$. Completing the square in the Gaussian
integral gives
\begin{equation}\label{eq:coarse-Gaussian-mgf}
 \frac1{a_N^{(t)}}\log\E\exp\{a_N^{(t)}\lambda\cdot p_\psi(\widehat X_N)\}
 =\sup_{u\in\cH_N}
       \{\lambda\cdot p_\psi(\widehat u_N)-q_N(u)\}.
\end{equation}
By \eqref{eq:coarse-sobolev}, $|\lambda\cdot p_\psi(\widehat u_N)|\le
C_{\lambda,\psi}\sqrt{q_N(u)}$, so near-maximizers have bounded $q_N$. The
compactness part of Lemma~\ref{lem:discrete-energy} bounds the upper limit of
the right side, and its recovery part, together with density of $C_c^\infty$
in $D^1$, bounds the lower limit. Thus the limit exists and equals
\begin{equation}\label{eq:coarse-limiting-mgf}
 \Lambda_\psi(\lambda)
 =\sup_{f\in D^1(\R^d)}
       \{\lambda\cdot p_\psi(f)-I(f)\},
\end{equation}
a finite quadratic form on $\R^k$. The G\"artner--Ellis theorem
\cite[Theorem~2.3.6(c)]{DemboZeitouni1998} applies, also when this form is
degenerate, and gives an LDP for $p_\psi(\widehat X_N)$ with rate
\begin{equation}\label{eq:coarse-projection-rate}
 J_\psi(z)=\inf\{I(f):f\in D^1(\R^d),\ p_\psi(f)=z\}.
\end{equation}
Indeed, the constraint defines a closed affine subspace of the Dirichlet
Hilbert space; its element of minimal norm identifies the Legendre transform
of \eqref{eq:coarse-limiting-mgf}, and in degenerate directions both sides
are infinite off the range of the covariance.

\smallskip
\noindent \emph{High modes on a fixed ball.} Fix $R$. Choose an orthonormal basis
$(e_j)_{j\ge1}$ of $L^2(B_R)$ with $e_j\in C_c^\infty(B_R)$, and let $P_k$
be the projection onto its first $k$ elements. Write
\[
 W_{N,k}=(I-P_k)(\widehat X_N|_{B_R}),\qquad
 \tau_{N,k}=\E\|W_{N,k}\|_2^2,
 \qquad \lambda_{N,k}=\|\Cov(W_{N,k})\|_{2\to2},
\]
where $\|\cdot\|_{2\to2}$ is the operator norm on $L^2(B_R)$. 
For fixed $N$, $k$, and $R$, the random element $W_{N,k}$ is a
centered Gaussian element of $L^2(B_R)$ with finite-rank covariance
operator (Indeed, $\widehat X_N|_{B_R}$ is piecewise constant on
finitely many lattice cubes, and applying $I-P_k$ preserves
finite-dimensionality); moreover, $\tau_{N,k}\to0$ for every
fixed $k$ by \eqref{eq:coarse-vanishing-trace} and $t_N\ge1$. Completing the
square for $\langle\widehat X_N,v\rangle$ and taking the supremum over unit
vectors $v$ orthogonal to the first $k$ modes gives the identity
\begin{equation}\label{eq:coarse-high-mode-duality}
 a_N^{(t)}\lambda_{N,k}
 =\frac12\sup_{q_N(u)\le1}
        \|(I-P_k)(\widehat u_N|_{B_R})\|_2^2.
\end{equation}
Set
\[
 \delta_k=\sup_{I(f)\le1}\|(I-P_k)(f|_{B_R})\|_2^2.
\]
The restrictions to $B_R$ of $\{I\le1\}$ form a compact subset of
$L^2(B_R)$, and $P_k\to I$ strongly, hence uniformly on this set; thus
$\delta_k\downarrow0$. Lemma~\ref{lem:discrete-energy} also implies
\begin{equation}\label{eq:coarse-high-mode-op}
 \limsup_N a_N^{(t)}\lambda_{N,k}\le\tfrac12\delta_k;
\end{equation}
otherwise, near-maximizers in \eqref{eq:coarse-high-mode-duality} would have
a subsequence converging in $L^2_{\rm loc}$ to a limit of energy at most one,
contradicting the definition of $\delta_k$.

Let $\gamma_k=\max\{\delta_k,k^{-1}\}>0$. For fixed $k$ and large $N$,
$a_N^{(t)}\lambda_{N,k}\le\gamma_k$. If $\rho_j$ are the eigenvalues of the
covariance of $W_{N,k}$, the Gaussian determinant formula and
$-\log(1-z)\le z/(1-z)$ give, for $2s\lambda_{N,k}<1$,
\[
 \log\E e^{s\|W_{N,k}\|_2^2}
 =-\tfrac12\sum_j\log(1-2s\rho_j)
 \le\frac{s\tau_{N,k}}{1-2s\lambda_{N,k}}.
\]
Taking $s=a_N^{(t)}/(4\gamma_k)$ and applying Markov's inequality yields
\[
 \limsup_N\frac1{a_N^{(t)}}\log\GP(\|W_{N,k}\|_2>\eps)
 \le-\frac{\eps^2}{4\gamma_k}.
\]
Since $\gamma_k\downarrow0$,
\begin{equation}\label{eq:coarse-high-mode}
 \lim_{k\to\infty}\limsup_{N\to\infty}\frac1{a_N^{(t)}}
  \log\GP\bigl(\|(I-P_k)(\widehat X_N|_{B_R})\|_2>\eps\bigr)
 =-\infty.
\end{equation}

\smallskip
\noindent \emph{The LDP on a fixed ball.} By \eqref{eq:coarse-projection-rate} and
the contraction principle, the laws of $P_k(\widehat X_N|_{B_R})$ satisfy an
LDP in $L^2(B_R)$ with rate
\[
 I_{R,k}(v)=\inf\{I(f):P_k(f|_{B_R})=v\}.
\]
By \eqref{eq:coarse-high-mode}, these are exponentially good approximations
of $\widehat X_N|_{B_R}$. The exponential approximation theorem
\cite[Theorem~4.2.16]{DemboZeitouni1998} therefore gives the LDP on
$L^2(B_R)$ with good rate
\begin{equation}\label{eq:coarse-ball-rate}
 I_R(v)=\inf\{I(f):f\in D^1(\R^d),\ f|_{B_R}=v\},
\end{equation}
once we check the identification of the rate and the closed-set condition of
that theorem. For each $M$, the set $\{f|_{B_R}:I(f)\le M\}$ is compact and
$\|(I-P_k)(f|_{B_R})\|_2\to0$ uniformly on it. Suppose $k_j\to\infty$ and
$I_{R,k_j}(v_j)$ stays bounded. Choose $f_j$ admissible for $I_{R,k_j}(v_j)$
within $1/j$ of the infimum, and pass to a subsequence along which the costs
converge to their lower limit. A further subsequence converges weakly in
$D^1$ and strongly in $L^2(B_R)$ to some $f$; the uniform projection estimate
gives $v_j\to f|_{B_R}$, and
$I_R(f|_{B_R})\le I(f)\le\liminf_jI_{R,k_j}(v_j)$. Conversely, if
$f|_{B_R}=v$, then $P_kv\to v$ and $I_{R,k}(P_kv)\le I(f)$. These two
observations identify the rate with $I_R$, show that it is good, and show
that $\inf_CI_R\le\liminf_k\inf_CI_{R,k}$ for every closed $C$.

\smallskip
\noindent 
\emph{Passage to $L^2_{\rm loc}$.} The space $L^2_{\mathrm{loc}}(\R^d)$, with
its usual topology, is the projective limit of the spaces $L^2(B_R)$,
$R=1,2,\ldots$, under restriction. The Dawson--G\"artner theorem
\cite[Theorem~4.6.1]{DemboZeitouni1998} gives an LDP there with rate
$\sup_RI_R(f|_{B_R})$, and this rate equals $I(f)$. The inequality
$\sup_RI_R(f|_{B_R})\le I(f)$ is immediate. Conversely, if the supremum is
$M<\infty$, choose $f_R\in D^1$ with $f_R=f$ on $B_R$ and
$I(f_R)\le M+1/R$; a weakly convergent subsequence in $D^1$ agrees with $f$
on every fixed ball, and lower semicontinuity gives $I(f)\le M$.
\end{proof}

\subsection{The bounded-potential Laplace principle}

\begin{proposition}\label{prop:laplace}
Let $V:\R\to\R$ be bounded and continuous. Then
\begin{equation}\label{eq:coarse-laplace}
 \lim_{N\to\infty}\frac1{a_N}
 \log\E\exp\left\{-\frac1{N^2\log N}
                   \sum_{x\in D_N}V(F_N(x))\right\}
 =-\inf_{f\in D^1(\R^d)}
      \left\{I(f)+\int_DV(f(y))\dd y\right\}.
\end{equation}
\end{proposition}

\begin{proof}
The map $\Psi^V(f)=\int_DV(f)$ is bounded and continuous on
$L^2_{\rm loc}$: convergence in $L^2(D)$ implies convergence in measure, and
then bounded continuity of $V$ gives convergence of $V(f)$ in $L^1(D)$. Let
$U_N=\bigcup_{x\in D_N}Q_{N,x}$. For every lattice function $u$,
\begin{equation}\label{eq:coarse-uniform-riemann}
 \left|N^{-d}\sum_{x\in D_N}V(u(x))-
                      \int_DV(\widehat u_N(y))\dd y\right|
 \le\|V\|_\infty|U_N\triangle D|=o(1),
\end{equation}
where the last step uses the Lipschitz boundary of $D$. Since
$(N^2\log N)^{-1}=a_NN^{-d}$, the exponent in \eqref{eq:coarse-laplace}
equals $-a_N(\Psi^V(\widehat F_N)+o(1))$, uniformly in the field. 

The claim
follows from Proposition~\ref{prop:coarse-ldp} and Varadhan's lemma for
bounded continuous functions \cite[Theorem~4.3.1]{DemboZeitouni1998}.
\end{proof}

\section{Proof in the critical window}\label{sec:critical-proof}

In this section we prove Theorem~\ref{thm:gff-window}. Thus
$\beta_N\to\beta\in\R$, and we write $h_N=h_N(\beta_N)$,
$\lambda_N=\theta e^{-\beta_N}$ and $\lambda=\theta e^{-\beta}$. We use the
skeleton $\sigma$-field $\cG_N$, the harmonic field $F_N=(h_N/g_0)m_N$, and
the sets $D_N^-$ and $R_N$ from Section~\ref{sec:coarse}.

\subsection{Upper bound}
To tackle the upper bound, we combine the conditional local estimate with the Laplace principle, and
then remove the truncation of the exponential potential.

\begin{proof}[Proof of the upper bound in Theorem~\ref{thm:gff-window}]
Fix $M$. Choose continuous functions $\chi_M:\R\to[0,1]$, increasing in $M$
to one, equal to one on $[-M,M]$ and to zero outside $[-M-1,M+1]$, and set
$V_M(t)=\lambda e^{-t}\chi_M(t)$.

\smallskip
\noindent \emph{Conditioning.} Discarding sites only increases the probability of the
hard-wall event. Conditionally on $\cG_N$, we keep in each box only the sites
of $D_N^-$ where $|F_N|\le M+1$; these sets are $\cG_N$-measurable, and
Lemma~\ref{lem:local-poisson} applies to them with $q_x=F_N(x)$. Each box
contributes a bounded intensity, at most $Ca_N$ boxes are involved, and
$\lambda_N\to\lambda$; by conditional independence and
$V_M(t)\le\lambda e^{-t}\mathbf1_{\{|t|\le M+1\}}$, we obtain
\begin{equation}\label{eq:critical-upper-local}
 \GP(\Omega_N(\beta_N))
 \le\E\exp\left\{-\frac1{N^2\log N}
                 \sum_{x\in D_N^-}V_M(F_N(x))+o(a_N)\right\},
\end{equation}
where the error is deterministic and uniform over the skeleton values. Since
$V_M$ is bounded and $|R_N|=o(N^d)$, extending the sum to all of $D_N$
changes the exponent by $o(a_N)$. Proposition~\ref{prop:laplace} therefore
gives
\begin{equation}\label{eq:critical-upper-var}
 \limsup_N a_N^{-1}\log\GP(\Omega_N(\beta_N))
 \le-\inf_{f\in D^1}J_M(f),\qquad
 J_M(f)=\alpha\cE_{\R^d}(f,f)+\int_DV_M(f)\dd x.
\end{equation}

\smallskip
\noindent \emph{Removing the truncation.} Since $V_M\le\lambda e^{-t}$ and $V_M$
increases with $M$, $\inf J_M$ is nondecreasing in $M$ and at most
$\cC_D^G(\beta)$. Choose $f_M$ with
$J_M(f_M)\le\inf J_M+M^{-1}\le\lambda|D|+M^{-1}$. The energies are bounded, so
along a subsequence $f_M$ converges weakly in $D^1$, in $L^2_{\rm loc}$, and
almost everywhere to some $f\in D^1$. At every point where $f$ is finite,
$|f_M|\le M$ eventually, so $V_M(f_M)\to\lambda e^{-f}$ almost everywhere.
Weak lower semicontinuity and Fatou's lemma give
\[
 \cC_D^G(\beta)
 \le\alpha\cE_{\R^d}(f,f)+\lambda\int_De^{-f}\dd x
 \le\liminf_M J_M(f_M),
\]
so $\inf J_M\uparrow\cC_D^G(\beta)$. Letting $M\to\infty$ in
\eqref{eq:critical-upper-var} proves the upper bound.
\end{proof}

\subsection{Lower bound}

We force the harmonic field into a tube around a prescribed profile, apply
the local estimate inside the tube, and restore the trimmed sites by
association.

\begin{proof}[Proof of the lower bound in Theorem~\ref{thm:gff-window}]
Fix a nonnegative $f\in C_c^\infty(\R^d)$ and $\eps>0$, and let
\[
 T_N(f,\eps)=
 \left\{\max_{x\in D_N^-}|F_N(x)-f(x/N)|\le\eps\right\}.
\]
By Lemma~\ref{lem:tube} with $t_N=1$,
$\GP(T_N(f,\eps))\ge\exp\{-a_N(\alpha\cE_{\R^d}(f,f)+o(1))\}$. On
$T_N(f,\eps)$ the values $F_N(x)$, $x\in D_N^-$, lie in
$[-\eps,\|f\|_\infty+\eps]$. Conditional independence,
Lemma~\ref{lem:local-poisson}, and a Riemann sum therefore give, on this
event,
\[
 \GP(\varphi_x\ge-h_N,\ x\in D_N^-\mid\cG_N)
 \ge\exp\left\{-a_N\left(e^\eps\lambda\int_De^{-f}\dd x+o(1)\right)\right\}.
\]
Integrating over the tube, we find that the hard wall on $D_N^-$ has
probability at least
$\exp\{-a_N(\alpha\cE_{\R^d}(f,f)+e^\eps\lambda\int_De^{-f}+o(1))\}$.

The covariances of the GFF are nonnegative, so by Gaussian association
\cite{Pitt1982},
\[
 \GP(\Omega_N(\beta_N))
 \ge\GP(\varphi_x\ge-h_N,\ x\in D_N^-)
       \prod_{x\in R_N}\GP(\varphi_x\ge-h_N).
\]
By Lemma~\ref{lem:one-point} and \eqref{eq:RN-size}, the logarithm of the
product is $-O(|R_N|/(N^2\log N))=-o(a_N)$. Letting $N\to\infty$ and then
$\eps\downarrow0$,
\[
 \liminf_N a_N^{-1}\log\GP(\Omega_N(\beta_N))
 \ge-\alpha\cE_{\R^d}(f,f)-\lambda\int_De^{-f}\dd x.
\]
Nonnegative functions in $C_c^\infty$ are dense in the nonnegative cone of
$D^1$ for the Dirichlet norm, hence also for local $L^2$ convergence, and
the exponential potentials then converge in $L^1(D)$ because they are
bounded by one. Taking the infimum over $f$ proves the lower bound in
\eqref{eq:critical-main}.

\smallskip
\noindent \emph{Local uniformity.} Restricting \eqref{eq:gff-variational} to $f\ge0$
shows that
$|\cC_D^G(\beta)-\cC_D^G(\beta')|\le\theta|D|\,|e^{-\beta}-e^{-\beta'}|$, so
the limit is continuous. Convergence along every sequence $\beta_N\to\beta$
to a continuous limit is equivalent to locally uniform convergence: if the
latter failed, there would be a compact set $K$, $\eps>0$, $N_k\to\infty$
and $\beta_k\in K$ with an error at least $\eps$ at $(N_k,\beta_k)$; after
extraction $\beta_k\to\beta_*$, and a sequence $\beta_N\to\beta_*$ through
the points $\beta_k$ would contradict \eqref{eq:critical-main}.
\end{proof}

\section{The off-window regimes}\label{sec:offwindow}

We now let the wall leave the window. The proofs follow the scheme of
Section~\ref{sec:critical-proof}, but every estimate has to be uniform at
the new, moving speed. We keep the skeleton, the boxes, the trimming radius
$r_N=(\log N)^{4/(d-2)}$, the harmonic field $F_N$, and the sets $D_N^-$ and
$R_N$ of Section~\ref{sec:coarse}, now with $h_N=h_N(\beta_N)$.

\subsection{The subcritical regime}\label{sec:off-subcritical}

Below the window, suppressing the violations forces the harmonic field above
a hard obstacle, and the minimal energy under this obstacle is the capacity
of $D$. Throughout this subsection,
\[
 t_N=-\beta_N\to\infty,\quad
 t_N^2=o(\log N),\quad h_N=h_N(\beta_N),\quad
 \Omega_N=\Omega_N(\beta_N),\quad G_N=\frac{F_N}{t_N}.
\]
Since $e^{t_N}=N^{o(1)}$ and $r_N/\ell_N=N^{-2/d+o(1)}$,
\eqref{eq:RN-size} gives
\begin{equation}\label{off:scales}
 e^{t_N}\frac{|R_N|}{N^2\log N}\le Ca_Ne^{t_N}\frac{r_N}{\ell_N}
 =o(a_Nt_N^2).
\end{equation}

\begin{lemma}[Moving one-point tails]\label{off:sub-tail}
For every fixed $C<\infty$, uniformly over boxes $B$, $x\in B^-$, and
$0\le q\le Ct_N$,
\begin{equation}\label{off:tail}
 \Pp\left(\eta_x^B<-h_N-\frac{g_0q}{h_N}\right)
 =\frac{\theta+o(1)}{N^2\log N}
   \exp\left\{t_N-q-\frac{g_0q^2}{2h_N^2}\right\}.
\end{equation}
Consequently, for fixed $\eps>0$, uniformly over skeleton values and over
sets $A\subset D_N^-$ on which $F_N\ge(1+\eps)t_N$,
\begin{equation}\label{off:favorable}
 \GP\left(\varphi_x\ge-h_N\text{ for }x\in A\mid\cG_N\right)
 \ge \exp\{-o(a_Nt_N^2)\}.
\end{equation}
\end{lemma}

\begin{proof}
For variance $g_0$, \eqref{off:tail} follows from Mills' ratio as in
Lemma~\ref{lem:one-point}, since $g_0q/h_N^2\to0$ in the stated range.
Replacing $g_0$ by $g_B(x,x)$ changes the exponent by
$O(h_N^2r_N^{2-d})=o(1)$, uniformly in this range.

For \eqref{off:favorable}, monotonicity in $q$ and \eqref{off:tail} bound
the conditional violation probability of each $x\in A$ by
$Ce^{-\eps t_N}/(N^2\log N)$. Conditionally on $\cG_N$, the residual fields
have nonnegative covariances (they are independent across boxes), so
Gaussian association bounds the conditional probability below by the product
of the one-site probabilities. Its logarithm is at least
$-Ca_Ne^{-\eps t_N}=-o(a_Nt_N^2)$.
\end{proof}

The next estimate forces the harmonic field above the obstacle. Its
uniformity in the skeleton values is essential when $t_N\to\infty$.

\begin{lemma}[Suppression of low regions]\label{off:low-regions}
Fix $0<\eps<1$ and $\delta>0$, and let
\[
 \cL_N(\eps,\delta)=
 \left\{\left|\{x\in D_N^-:F_N(x)\le(1-\eps)t_N\}\right|
                         \ge\delta N^d\right\}.
\]
There is a deterministic sequence $w_N$ with $w_N/t_N^2\to\infty$ such that
\begin{equation}\label{off:low-regions-bound}
 \GP(\Omega_N\cap\cL_N(\eps,\delta))\le e^{-a_Nw_N}.
\end{equation}
\end{lemma}

\begin{proof}
Put $q_0=(1-\eps)t_N$, and for a box $B$ let
$S_B=\{x\in D_N\cap B^-:F_N(x)\le q_0\}$. At most $Ca_N$ boxes meet $D_N$,
each with at most $\ell_N^d$ sites. Hence, on $\cL_N(\eps,\delta)$, at least
$c_\delta a_N$ boxes satisfy $|S_B|\ge c_\delta\ell_N^d$, for some
$c_\delta>0$; we call these boxes \emph{rich}. The rich boxes and their sets
$S_B$ are $\cG_N$-measurable.

Condition on $\cG_N$ and fix a rich box $B$. On $S_B$ we have $F_N\le q_0$,
so the wall constraint there implies the weaker constraint with $F_N$ replaced
by $q_0$. It therefore suffices to bound
\begin{equation}\label{off:one-box}
 \Pp\left(\eta_x^B\ge-h_N-\frac{g_0q_0}{h_N},\ x\in S_B\right).
\end{equation}
By \eqref{off:tail} and $t_N^2=o(\log N)$, the sum of the one-site
violation probabilities over $S_B$ is at least
\[
 \mu_{S_B}\ge
 \frac{c|S_B|}{N^2\log N}e^{\eps t_N}
 \ge c_\delta e^{\eps t_N}.
\]
The standardized thresholds $z_x=(h_N+g_0q_0/h_N)/\sqrt{g_B(x,x)}$ are
positive and satisfy $z_x^2=4\log N+o(\log N)$ uniformly, so
Lemma~\ref{lem:mesoscopic-comparison} bounds \eqref{off:one-box} by
\[
 e^{-c_\delta e^{\eps t_N}}+N^{-c+o(1)}\le e^{-v_N},\qquad
 v_N\asymp\min\{e^{\eps t_N},\log N\},\qquad
 \frac{v_N}{t_N^2}\longrightarrow\infty .
\]
By conditional independence of the boxes,
$\GP(\Omega_N\mid\cG_N)\le e^{-c_\delta a_Nv_N}$ on $\cL_N(\eps,\delta)$.
Taking expectations proves \eqref{off:low-regions-bound} with
$w_N=c_\delta v_N$.
\end{proof}

\begin{proof}[Proof of Theorem~\ref{thm:subcritical-sharp}]
Since $h_N^2\sim4g_0\log N$, Proposition~\ref{prop:coarse-ldp} gives an
LDP for $\widehat G_N=\widehat F_N/t_N$ in $L^2_{\rm loc}(\R^d)$, at speed
$a_Nt_N^2$, with good rate function $I$.

\smallskip
\noindent \emph{Upper bound.} Fix $0<\eps<1/4$ and apply
Lemma~\ref{off:low-regions} with $\delta=\eps$. Up to an event that is
negligible at speed $a_Nt_N^2$, the hard wall implies $G_N>1-\eps$ at all but
$\eps N^d$ sites of $D_N^-$. By \eqref{eq:RN-size} and
$|U_N\triangle D|\to0$ (with $U_N$ as in the proof of
Proposition~\ref{prop:laplace}), for large $N$ the interpolation
$\widehat G_N$ then belongs to
\[
 A_\eps=\{f\in L^2_{\rm loc}(\R^d):
                 |\{y\in D:f(y)\le1-\eps\}|\le2\eps\}.
\]
If $f_n\in A_\eps$ and $f_n\to g$ in $L^2(D)$, then
$\{g<1-2\eps\}\subset\{f_n\le1-\eps\}\cup\{|f_n-g|>\eps\}$, and the last set
has vanishing measure; hence the closure of $A_\eps$ is contained in
$\{g:|\{y\in D:g(y)<1-2\eps\}|\le2\eps\}$. The LDP upper bound for closed
sets yields
\begin{equation}\label{off:sub-upper}
 \limsup_{N\to\infty}\frac{\log\GP(\Omega_N)}{a_Nt_N^2}
 \le-\inf\{I(g):|\{y\in D:g(y)<1-2\eps\}|\le2\eps\}.
\end{equation}
As $\eps\downarrow0$, the infimum on the right tends to $\alpha\capR(D)$. It
is clearly at most $\alpha\capR(D)$. Conversely, take approximate minimizers
$g_\eps$ along a sequence $\eps\downarrow0$. They have bounded energy, so by
Sobolev and Rellich compactness a subsequence converges in $L^2_{\rm loc}$
and almost everywhere to some $g$. For fixed $\eta>0$ and $\eps<\eta$, the
set $\{g_\eps<1-2\eta\}\cap D$ has measure at most $2\eps$, so by Fatou's
lemma $g\ge1-2\eta$ a.e.\ on $D$. As $\eta$ is arbitrary, $g\ge1$ a.e.\ on
$D$, and lower semicontinuity of the energy gives
$\alpha\capR(D)\le I(g)\le\liminf I(g_\eps)$.

\smallskip
\noindent \emph{Lower bound.} Fix $\eps>0$ and $u\in C_c^\infty(\R^d)$ with
$u\ge1+3\eps$ near $\overline D$. By Lemma~\ref{lem:tube},
\[
 \GP\left(\max_{x\in D_N^-}|G_N(x)-u(x/N)|\le\eps\right)
 \ge\exp\left\{-a_Nt_N^2
       \left(\alpha\,{\cE_{\R^d}(u,u)}+o(1)\right)\right\}.
\]
On this $\cG_N$-measurable event, $F_N\ge(1+2\eps)t_N$ on $D_N^-$, so
\eqref{off:favorable} gives
\[
 \GP(\varphi_x\ge-h_N\text{ for }x\in D_N^-)
 \ge\exp\left\{-a_Nt_N^2
       \left(\alpha\,{\cE_{\R^d}(u,u)}+o(1)\right)\right\}.
\]
To restore $R_N$, note that $\GP(\varphi_x<-h_N)\le Ce^{t_N}/(N^2\log N)$
by \eqref{off:tail}. By Gaussian association and \eqref{off:scales}, the wall
on $R_N$ has probability at least
\[
 \exp\left\{-C\frac{e^{t_N}|R_N|}{N^2\log N}\right\}
 =\exp\{-o(a_Nt_N^2)\},
\]
and a further application of association combines the two increasing wall
events. Optimizing over $u=(1+3\eps)h$ with $h$ as in
Lemma~\ref{lem:smooth-capacity} gives the
lower bound $-(1+3\eps)^2\alpha\capR(D)$, and letting $\eps\downarrow0$
completes the proof.
\end{proof}

\begin{remark}\label{rem:sub-range}
The restriction $t_N^2=o(\log N)$ enters only through
Lemma~\ref{off:low-regions}: the additive Gaussian comparison error in a rich
box is $N^{-c+o(1)}=e^{-c'\log N}$, and it must be negligible compared with
the one-box cost scale $t_N^2$. Extending the theorem to the full range
$t_N=o(\log N)$ would require a local estimate that remains sharp when the
local intensity grows.
\end{remark}

\subsection{The supercritical regime}\label{sec:off-supercritical}

Above the window, the one-point violation intensity gives the leading cost.
Throughout this subsection $\beta_N\to+\infty$ and $\beta_N=o(\log N)$, and
\[
 h_N=h_N(\beta_N),\qquad \lambda_N=\theta e^{-\beta_N},\qquad
 a_N^{+}=a_N\lambda_N,\qquad \Omega_N=\Omega_N(\beta_N).
\]
Then $\lambda_N\to0$ but $\lambda_N=N^{-o(1)}$, so $a_N^{+}\to\infty$. All
errors below are measured at the speed $a_N^{+}$, which is smaller than the
critical speed $a_N$.

\begin{lemma}[Small-intensity local estimate]\label{off:super-local}
For fixed $M<\infty$, uniformly over boxes $B$, sets $A\subset B^-$, and
deterministic $q:A\to[-M,M]$,
\begin{equation}\label{off:super-poisson}
 \log\Pp\left(\eta_x^B\ge-h_N-\frac{g_0q_x}{h_N}
                                      \text{ for }x\in A\right)
 =-\frac{\lambda_N}{N^2\log N}\sum_{x\in A}e^{-q_x}
      +o(\lambda_N).
\end{equation}
\end{lemma}

\begin{proof}
Mills' ratio \eqref{eq:mills} and $h_N^2(g_0-g_B(x,x))=o(1)$ give, uniformly,
\[
 p_x:=\Pp\left(\eta_x^B<-h_N-\frac{g_0q_x}{h_N}\right)
 =(1+o(1))\frac{\lambda_Ne^{-q_x}}{N^2\log N},
 \qquad \sum_{x\in A}p_x=O(\lambda_N).
\]
The squared standardized thresholds are
$4\log N+\log\log N+2\beta_N+O(1)=4\log N+o(\log N)$, so
Lemma~\ref{lem:mesoscopic-comparison} gives an additive comparison error
$N^{-c+o(1)}=o(\lambda_N)$. The probability for independent sites equals
$\exp\{-\sum p_x+O(\sum p_x^2)\}$, and $\sum p_x^2=o(\lambda_N)$. Since these
probabilities tend to one uniformly, taking logarithms preserves the
$o(\lambda_N)$ error.
\end{proof}

\begin{proof}[Proof of Theorem~\ref{thm:supercritical-sharp}]
\emph{Lower bound.} The one-point Mills expansion is
$p_N:=\GP(\varphi_0<-h_N)=(1+o(1))\lambda_N/(N^2\log N)$. Gaussian association
gives directly
\begin{equation}\label{off:super-lower}
 \log\GP(\Omega_N)\ge |D_N|\log(1-p_N)
 =-(1+o(1))|D|a_N^{+}.
\end{equation}

\smallskip
\noindent \emph{Upper bound.} Fix $M>0$ and $R>M$, and introduce the bounded
Lipschitz function and the empirical potential
\[
 W_M(s)=e^{-\max\{s,-M\}},\qquad
 \cV_{N,M}=\frac1{N^d}\sum_{x\in D_N^-}W_M(F_N(x)).
\]
Conditionally on $\cG_N$, raise every value of $F_N$ below $-M$ to $-M$, and
keep only the sites of $D_N^-$ where the resulting value is at most $R$. Both
operations increase the conditional probability of the wall. Applying
Lemma~\ref{off:super-local} (with $M$ replaced by $R$) in each box and using
conditional independence gives, uniformly in the skeleton values,
\begin{equation}\label{off:super-conditional}
 \GP(\Omega_N\mid\cG_N)
 \le\exp\{-a_N^{+}\cV_{N,M}+o_R(a_N^{+})+Ca_N^{+}e^{-R}\}.
\end{equation}
Here $o_R(a_N^+)$ collects the local errors of the at most $Ca_N$ relevant
boxes, and the last term bounds the contribution to $\cV_{N,M}$ of the
discarded sites, where $F_N>R$.

It remains to show that $\cV_{N,M}$ concentrates around $|D|$ at speed
$a_N\gg a_N^+$. By \eqref{eq:coarse-uniform-riemann} and
\eqref{eq:RN-size},
\[
 \cV_{N,M}=\int_D W_M(\widehat F_N(y))\dd y+o(1),
\]
uniformly in the field. The integral is a continuous function on
$L^2_{\rm loc}$ which equals $|D|$ at the zero function, and the good rate
function $I$ of Proposition~\ref{prop:coarse-ldp} vanishes only there. By
compactness of its sublevel sets, $I$ is bounded below by some $2c_\delta>0$
on the closed set where the integral differs from $|D|$ by at least
$\delta/2$. Hence, for every fixed $\delta>0$ and all large $N$,
\[
 \GP\bigl(|\cV_{N,M}-|D||>\delta\bigr)\le e^{-c_\delta a_N}.
\]
Splitting the expectation of \eqref{off:super-conditional} according to
this event, and using $\cV_{N,M}\ge0$,
\[
 \E e^{-a_N^{+}\cV_{N,M}}
 \le e^{-a_N^{+}(|D|-\delta)}+e^{-c_\delta a_N}.
\]
Since $a_N^{+}/a_N=\lambda_N\to0$, the second term is negligible at speed
$a_N^{+}$. Letting $N\to\infty$, then $R\to\infty$, then $\delta\downarrow0$,
we obtain
$\limsup_{N}(a_N^{+})^{-1}\log\GP(\Omega_N)\le-|D|$. Together with
\eqref{off:super-lower} this proves \eqref{eq:supercritical-main}.
\end{proof}

\section{Interlacement cover levels}\label{sec:ri}

We now apply the hard-wall estimates to cover levels of random
interlacements. The upper bounds come from a comparison with the GFF; the
lower bounds use an entropy identity and a local capacity estimate, followed
by sprinkling.

For a finite set $K\subset\Z^d$, we have
\[
 \mathbb P(K\cap\mathcal I^u=\varnothing)
   =\exp\{-u\capZ(K)\},
 \qquad\text{in particular}\qquad
 \mathbb P(x\notin\mathcal I^u)=e^{-u/g_0},
\]
with $\capZ$ as in Section~\ref{sec:normalizations}. We set
\[
 v_N=dg_0\log N,
 \qquad
 \Psi_N=a_N(\log\log N)^2.
\]

\subsection{The change-of-measure input and a local estimate}
To state the entropy identity, let $\mathbb P^u$ denote the standard interlacement
law at level $u$ with independent mean-one exponential holding times
attached to the trajectories; its trace has the usual interlacement law. For
a positive function $f$ on $\Z^d$ equal to $1$ outside a finite set, let
$\operatorname{cap}_f$ denote the electrical capacity for the conductances
$c_f(x,y)=f(x)f(y)/(2d)$ on nearest-neighbour edges.

\begin{lemma}[Change of measure]\label{lem:ri-input}
For $u>0$ and such a function $f$, there is a probability law
$\widetilde{\mathbb P}^{u,f}\ll\mathbb P^u$ on interlacement trajectories
whose trace $\mathcal I$ satisfies
\begin{align}
 H(\widetilde{\mathbb P}^{u,f}\mid\mathbb P^u)
   &=u\,\cE_{\mathbb Z^d}(f-1,f-1),\label{eq:ri-entropy}\\
 \widetilde{\mathbb P}^{u,f}(K\cap\mathcal I=\varnothing)
   &=\exp\{-u\operatorname{cap}_f(K)\}
       \qquad(K\subset\mathbb Z^d\text{ finite}).\label{eq:ri-void}
\end{align}
\end{lemma}

\begin{proof}
These identities follow from the change-of-measure computations in
\cite[Lemma~1.2, Proposition~2.1, (2.9), and the proof of
Proposition~2.3]{LiSznitman2014}, which apply without change to every
positive finite-support perturbation of $1$; we refer to that article for
the construction of the tilted law, and recall only what is used. Under
$\widetilde{\mathbb P}^{u,f}$ the trajectories are those of the
continuous-time walk with jump rates $f(y)/(2df(x))$ from $x$ to a neighbour
$y$. This walk is reversible with respect to $f^2$, its conductances are
$f^2(x)\cdot f(y)/(2df(x))=c_f(x,y)$, and the tilted interlacements have
occupation-time intensity $uf(x)^2$ at $x$:
\[
 \widetilde{\mathbb E}^{u,f}\Bigl[\,\text{total time spent at }x\,\Bigr]
 =u\,f(x)^2 .
\]
The void probability \eqref{eq:ri-void} is the usual one for interlacements
of a reversible walk, with $\operatorname{cap}_f(K)=\sum_{x\in K}\widetilde
e_K(x)$ the total mass of the equilibrium measure of $K$ for this walk
\cite[(1.49)--(1.50)]{LiSznitman2014}; by the variational characterization of
capacity, this total mass is the electrical capacity of $K$ for the
conductances $c_f$. With $\Delta_{\rm dis}=P-I$, the entropy computation of
\cite{LiSznitman2014} reads $u\sum_x f(x)(-\Delta_{\rm dis}f)(x)$. Writing
$f=1+\phi$ with $\phi$ finitely supported, $\sum_x(-\Delta_{\rm
dis}\phi)(x)=0$, so this equals
$u\langle\phi,(I-P)\phi\rangle=u\,\cE_{\Z^d}(f-1,f-1)$.
\end{proof}

\begin{lemma}\label{lem:ri-singleton}
Let $f\ge1$, $1\le s\le S$, and suppose that $f=s$ on the lattice ball of
radius $r\ge2$ about $x$. Then, with $C$ depending only on $d$ and $S$,
\begin{equation}\label{eq:ri-singleton}
 \operatorname{cap}_f(\{x\})\geq \frac{s^2}{g_0}-Cr^{2-d}.
\end{equation}
Consequently, if $f=s$ on a fixed macroscopic neighbourhood of $D_N$ and
$u=O(\log N)$, then uniformly in $x\in D_N$,
\begin{equation}\label{eq:ri-local-vacancy}
 \widetilde{\mathbb P}^{u,f}(x\notin\mathcal I)
 \leq\exp\{-us^2/g_0+o(1)\}.
\end{equation}
\end{lemma}

\begin{proof}
We bound the effective resistance from $x$ to infinity by Thomson's
principle (see e.g.\ \cite[Section~2.4]{LyonsPeres2016}). For the original conductances $c_0(y,z)=1/(2d)$, the unit current
from $x$ to infinity is
\[
 j(y,z)=c_0(y,z)\bigl(g(x,y)-g(x,z)\bigr),
\]
and its energy $\sum_{\rm edges}j^2/c_0$ equals $\cE_{\Z^d}(g(x,\cdot),g(x,\cdot))=g_0$.
Let $B$ be the ball of radius $r/2$ about $x$, extend $g_B(x,\cdot)$ by
zero, and put $w=g(x,\cdot)-g_B(x,\cdot)$. Since $w$ is harmonic in $B$,
Dirichlet orthogonality gives
\[
 \cE_{\mathbb Z^d}(w,w)=g_0-g_B(x,x)\leq Cr^{2-d}.
\]
For $r\geq4$, every edge not contained in the radius-$r$ plateau has both
endpoints outside $B$, where the gradients of $g(x,\cdot)$ and $w$ agree.
Hence the current $j$ has $c_0$-energy at most $Cr^{2-d}$ off the plateau.
On the plateau, $c_f=s^2c_0$, while $c_f\ge c_0$ everywhere because $f\ge1$.
Thomson's principle gives
\[
 \operatorname{cap}_f(\{x\})^{-1}
  \leq \frac{g_0}{s^2}+Cr^{2-d}.
\]
Inverting, uniformly for $1\leq s\leq S$, and enlarging $C$ to cover
$2\leq r<4$, proves \eqref{eq:ri-singleton}. Combined with \eqref{eq:ri-void}
and $(\log N)N^{2-d}\to0$, it gives \eqref{eq:ri-local-vacancy}.
\end{proof}

\subsection{Lower bounds by tilting}

\begin{proof}[Proof of the lower bounds in Proposition~\ref{prop:ri-phase}]
Choose $h\in C_c^\infty(\mathbb R^d)$ with $0\leq h\leq1$ and $h=1$
near $D$, and write $h^{[N]}(x)=h(x/N)$. For levels $0<u_1<v$, apply
Lemma~\ref{lem:ri-input} with
\[
 f_N=1+\left(\sqrt{v/u_1}-1\right)h^{[N]},
\]
and write $\widetilde{\mathbb P}$, $\widetilde{\mathbb E}$ for
$\widetilde{\mathbb P}^{u_1,f_N}$ and its expectation. Since
$f_N=\sqrt{v/u_1}$ near $D_N$, the tilted interlacements look locally like
interlacements at level $v$. Precisely, if $v/u_1$ stays bounded and
$u_1=O(\log N)$, Lemmas~\ref{lem:ri-input} and~\ref{lem:ri-singleton} give
\begin{equation}\label{eq:ri-holes}
 \begin{aligned}
 \widetilde{\mathbb E}|D_N\setminus\mathcal I|
       &\leq C N^d e^{-v/g_0},\quad
 \mathsf H_N:=H(\widetilde{\mathbb P}^{u_1,f_N}\mid\mathbb P^{u_1})
   &=(\sqrt v-\sqrt{u_1})^2
             \cE_{\mathbb Z^d}(h^{[N]},h^{[N]}).
 \end{aligned}
\end{equation}
By
\eqref{eq:energy-scaling},
\begin{equation}\label{eq:ri-energy-scale}
 \cE_{\mathbb Z^d}(h^{[N]},h^{[N]})
   =\frac{N^{d-2}}{d}\,\cE_{\R^d}(h,h)+o(N^{d-2}).
\end{equation}

Let $\mathcal A_N=\{|D_N\setminus\mathcal I|\leq k_N\}$ for a sequence $k_N$
chosen below, and suppose that $\widetilde{\mathbb P}(\mathcal A_N)\to1$.
The entropy inequality gives
\begin{equation}\label{eq:ri-event-entropy}
 \log\mathbb P^{u_1}(\mathcal A_N)
 \geq-\frac{\mathsf H_N+\log2}{\widetilde{\mathbb P}(\mathcal A_N)}.
\end{equation}
We add independent interlacements at level $u_2>0$. For a fixed set $S$ of at
most $k_N$ vertices, the events $\{x\in\mathcal I^{u_2}\}$, $x\in S$, are
increasing functions of the Poisson cloud of trajectories, and only the
finitely many trajectories that hit $S$ matter. By the Harris--FKG
inequality for Poisson processes \cite[Section~20.3]{LastPenrose2017}, $S$ is
covered by the new trajectories with probability at least
$\prod_{x\in S}\Pp(x\in\mathcal I^{u_2})\ge(1-e^{-u_2/g_0})^{k_N}$.
Conditioning on the first trace, we may apply this to its vacant set. Since the union of the two independent
traces is an interlacement at level $u_1+u_2$,
\begin{equation}\label{eq:ri-sprinkling}
 \mathbb P(\mathcal T_N\leq u_1+u_2)
 \geq \mathbb P^{u_1}(\mathcal A_N)(1-e^{-u_2/g_0})^{k_N}.
\end{equation}

\noindent \emph{Subcritical levels.} Fix $\gamma\in(0,2/d)$ and
$0<\varepsilon<\min\{d\gamma,d-2\}$, and put
\[
 L=\log N,\qquad
 u_1=g_0(d\gamma-\varepsilon)L,\quad
 u_2=g_0\varepsilon L,\quad
 v=g_0(2+\varepsilon)L,\quad
 k_N=\left\lceil N^{d-2-\varepsilon/2}\right\rceil.
\]
Then $u_1+u_2=\gamma v_N$. By Markov's inequality and \eqref{eq:ri-holes},
$\widetilde{\mathbb P}(\mathcal A_N^c)\le CN^{-\varepsilon/2}\to0$. The
logarithm of the second term in the right-hand side of \eqref{eq:ri-sprinkling} is
$O(N^{d-2-3\varepsilon/2})=o(N^{d-2}L)$. By \eqref{eq:ri-event-entropy} and
\eqref{eq:ri-energy-scale},
\[
 \liminf_N\frac{\log\mathbb P(\mathcal T_N\leq\gamma v_N)}{N^{d-2}L}
 \geq-{\frac{g_0}{d}}
   \left(\sqrt{2+\varepsilon}-\sqrt{d\gamma-\varepsilon}\right)^2
     {\cE_{\R^d}(h,h)}.
\]
Letting $\varepsilon\downarrow0$, using
${\frac{1}{d}(\sqrt2-\sqrt{d\gamma})^2=(\sqrt{2/d}-\sqrt\gamma)^2}$,
and taking the infimum over $h$ proves the lower bound in
\eqref{eq:ri-subcritical}.

\noindent \emph{The transition level.} Set $L=\log N$, $\ell=\log L=\log\log N$, and
\[
 s_N=\tfrac12\ell,\qquad
 u_1=g_0(2L-s_N),\quad u_2=g_0s_N,\quad v=g_0(2L+s_N),\qquad
 k_N=\left\lceil\ell N^{d-2}e^{-s_N}\right\rceil.
\]
Then $u_1+u_2=2g_0L$, and \eqref{eq:ri-holes} with Markov's inequality gives
$\widetilde{\mathbb P}(\mathcal A_N)\geq1-C/\ell$. Moreover,
\[
 (\sqrt v-\sqrt{u_1})^2
  =\frac{g_0s_N^2}{2L}(1+o(1)),
 \qquad\text{hence}\qquad
 \lim_N\frac{\mathsf H_N}{\Psi_N}
  ={\frac{g_0}{8d}\,\cE_{\R^d}(h,h)},
\]
where $\Psi_N=a_N(\log\log N)^2=\frac{N^{d-2}}{\log N}(\log \log N)^2$.
Since $e^{-s_N}=L^{-1/2}$, the logarithm of the second term in the right-hand side of \eqref{eq:ri-sprinkling} satisfies
\[
 -k_N\log(1-e^{-s_N})
  \leq 2k_Ne^{-s_N}
  =O\left(\frac{\ell N^{d-2}}{L}\right)=o(\Psi_N).
\]
Applying \eqref{eq:ri-event-entropy}--\eqref{eq:ri-sprinkling} and
optimizing over $h$ gives \eqref{eq:ri-critical-lower}.
\end{proof}

\subsection{Upper bounds and the predicted window}

\begin{proof}[Proof of the upper bounds in Proposition~\ref{prop:ri-phase}]
Lupu's coupling \cite[Theorem~3]{Lupu2016} gives
$\{x:\varphi_x>\sqrt{2u}\}\subset\mathbb Z^d\setminus\mathcal I^u$.
Hence, for every finite $A$, using the symmetry of $\varphi$,
\begin{equation}\label{eq:ri-gff-comparison}
 \mathbb P(A\subset\mathcal I^u)
 \leq\GP\left(\varphi_x\leq\sqrt{2u}\text{ for all }x\in A\right)
 =\GP\left(\varphi_x\geq-\sqrt{2u}\text{ for all }x\in A\right).
\end{equation}
Let $D^\circ$ be the interior of $D$; then $D^\circ_N=(ND^\circ)\cap\Z^d\subset D_N$, so
the hard wall on $D_N$ implies the hard wall on $D^\circ_N$.

\smallskip
\noindent \emph{Subcritical levels.} Write
$p_N(a)=\GP(\varphi_x\ge-a\sqrt{\log N},\ x\in D^\circ_N)$. Shifting the
field by $a\sqrt{\log N}\,h^{[N]}$, which equals $a\sqrt{\log N}$ on
$D^\circ_N$, and applying the Cameron--Martin formula and H\"older's
inequality with exponents $s$ and $s/(s-1)$, we get for $s>1$
\[
 p_N(a)\le p_N(0)^{1/s}
 \exp\left\{\frac{a^2\log N}{2(s-1)}
                   \cE_{\mathbb Z^d}(h^{[N]},h^{[N]})\right\}.
\]
In our normalization, the zero-wall asymptotic (see \cite[Theorem~1 and Section~1.6.3]{FeldheimFeldheimMuirhead2026};
see also \cite{BDZ1995}), reads
\[
 \log p_N(0)=-(2g_0+o(1))\capZ(D^\circ_N)\log N.
\]
By \eqref{eq:capacity-scaling} and $\capR(D^\circ)=\capR(D)$, the zero-wall
rate at speed $N^{d-2}\log N$ is $C_0={2g_0\capR(D)/d}$. Combining this with
\eqref{eq:ri-energy-scale} and minimizing over $h$ yields
\[
 \limsup_N\frac{\log p_N(a)}{N^{d-2}\log N}
 \le -\frac{C_0}{s}+\frac{A}{s-1},\qquad
 C_0={\frac{2g_0}{d}}\capR(D),\quad A={\frac{a^2}{2d}}\capR(D).
\]
For $0<a<2\sqrt{g_0}$ we have $A<C_0$, and optimizing over $s>1$ gives
$-(\sqrt{C_0}-\sqrt A)^2=-(2\sqrt{g_0}-a)^2\capR(D)/(2d)$. With
$a=\sqrt{2dg_0\gamma}$, so that $a\sqrt{\log N}=\sqrt{2\gamma v_N}$,
\eqref{eq:ri-gff-comparison} gives the upper bound in
\eqref{eq:ri-subcritical}.

\smallskip
\noindent \emph{The transition level.} For $u=2g_0\log N$, the level $\sqrt{2u}$ in
\eqref{eq:ri-gff-comparison} equals $h_N(\beta_N)$ with
$\beta_N=-\tfrac12\log\log N$, which satisfies $\beta_N^2=o(\log N)$.
Theorem~\ref{thm:subcritical-sharp}, applied to 
$D^\circ$, gives \eqref{eq:ri-critical-upper}: the factor $\beta_N^2$
contributes $\frac14(\log\log N)^2$, which turns ${g_0/(8d)}$ into
${g_0/(32d)}$. This proves Proposition~\ref{prop:ri-phase}.
\end{proof}

\begin{remark}[The RI window and a phase-transition prediction]\label{rem:ri-window}
The one-point vacancy probability $e^{-u/g_0}$ carries no Gaussian
prefactor. Consequently, at
\begin{equation}\label{eq:ri-predicted-window}
 u_N^{\rm RI}(\beta)=g_0(2\log N+\log\log N+\beta),
\end{equation}
the expected number of vacant vertices in $D_N$ is asymptotically
$|D|e^{-\beta}{a_N}$. By contrast, under $u=h^2/2$ the GFF window
\eqref{eq:window} corresponds to
\[
 g_0(2\log N+\tfrac12\log\log N+\beta).
\]
The different coefficient of $\log\log N$ reflects the Gaussian Mills
prefactor, and it is the reason why the comparison
\eqref{eq:ri-gff-comparison} does not identify the RI window rate.

The off-critical change of mechanism was already described in
\cite[Section~0.1]{LiShiXu2024}; the finer centering and the rate proposed
below are not asserted there. We expect that, for each fixed $\beta\in\R$,
\begin{equation}\label{eq:ri-window-prediction}
 \lim_{N\to\infty}{\frac1{a_N}}
    \log\Pp(\mathcal T_N\le u_N^{\rm RI}(\beta))
 =-\cC_D^{\rm RI}(\beta),
\end{equation}
with the predicted rate
\begin{equation}\label{eq:ri-rate-prediction}
 \cC_D^{\rm RI}(\beta)=\inf_{F\in D^1(\R^d)}
 \left\{{\frac{g_0}{8d}\,\cE_{\R^d}(F,F)}
          +e^{-\beta}\int_D e^{-F(x)}\dd x\right\}.
\end{equation}
Formally, raising the local intensity by $g_0F$ multiplies the vacancy
intensity by $e^{-F}$. Expanding the entropy \eqref{eq:ri-entropy} to second
order at intensity $2g_0\log N$, with $f^2=1+F/(2\log N)$, yields the
{coefficient $g_0/(8d)$ in front of $\cE_{\R^d}(F,F)$ at speed $a_N$}. The competition
between a collective profile and the remaining vacancies would thus be the
same as for the GFF, with a different centering and a different local
coefficient. A proof would have to control vacancy correlations at this
speed.

Matching this prediction with the subcritical regime also suggests
\begin{equation}\label{eq:ri-critical-prediction}
 \lim_{N\to\infty}
 \frac{\log\Pp(\mathcal T_N\le2g_0\log N)}
 {{a_N(\log\log N)^2}}
 =-{\frac{g_0}{8d}}\capR(D),
\end{equation}
which would make the lower bound \eqref{eq:ri-critical-lower} sharp: the
level $2g_0\log N$ corresponds to $\beta=-\log\log N$ in
\eqref{eq:ri-predicted-window}, twice the GFF value. 
Both predictions \eqref{eq:ri-window-prediction} and \eqref{eq:ri-critical-prediction} remain open.
\end{remark}

\phantomsection

\end{document}